\documentclass{svjour3}

\usepackage{amsmath,amsfonts,amssymb}
\usepackage{graphicx}
\usepackage{bm}
\usepackage{hyperref}
\usepackage{rotating}
\usepackage{float}
\usepackage{color}
\usepackage{soul}

\usepackage{subfigure}

\setkeys{Gin}{draft=false}
\usepackage{natbib}

\authorrunning{DAVIS and LEVEQUE}

\titlerunning{Adjoint Methods in Tsunami Modeling}

\journalname{Submitted}

\graphicspath{{../figs/}}

\usepackage{color}
\definecolor{red}{rgb}{1.0,0.,0.}

\definecolor{darkblue}{rgb}{0.2,0.2,1.0}
\newcommand{\revised}[2]{#2}  
\newcommand{\removed}[1]{}
\newcommand{\added}[1]{#1}
\definecolor{darkgreen}{rgb}{0.1,0.5,0.1}

\newcommand\reals{{{\rm l} \kern -.15em {\rm R} }}
\newcommand{\ignore}[1]{}

\newcommand{\cref}[1]{(\ref{#1})}
\newcommand{\Cref}[1]{(\ref{#1})}
\newcommand{\Fig}[1]{Figure~\ref{#1}}
\newcommand{\Sec}[1]{Section~\ref{#1}}

\ignore{
\newcommand*\patchAmsMathEnvironmentForLineno[1]{%
  \expandafter\let\csname old#1\expandafter\endcsname\csname #1\endcsname
  \expandafter\let\csname oldend#1\expandafter\endcsname\csname end#1\endcsname
  \renewenvironment{#1}%
     {\linenomath\csname old#1\endcsname}%
     {\csname oldend#1\endcsname\endlinenomath}}%
\newcommand*\patchBothAmsMathEnvironmentsForLineno[1]{%
  \patchAmsMathEnvironmentForLineno{#1}%
  \patchAmsMathEnvironmentForLineno{#1*}}%
\AtBeginDocument{%
\patchBothAmsMathEnvironmentsForLineno{equation}%
\patchBothAmsMathEnvironmentsForLineno{align}%
\patchBothAmsMathEnvironmentsForLineno{flalign}%
\patchBothAmsMathEnvironmentsForLineno{alignat}%
\patchBothAmsMathEnvironmentsForLineno{gather}%
\patchBothAmsMathEnvironmentsForLineno{multline}%
}
}

\begin{document}

\title{Adjoint Methods for Guiding Adaptive Mesh Refinement in Tsunami
Modeling%
  \thanks{Supported in part by an NSF Graduate Research Fellowship 
  DGE-1256082 and NSF grants DMS-1216732 and EAR-133141.}}

\author{B. N. Davis and R. J. LeVeque}

\institute{%
Department of Applied Mathematics, University of Washington,
Seattle, WA. 
}

\maketitle


\begin{abstract}
One difficulty in developing numerical methods for
tsunami modeling is the fact 
that solutions contain \added{time-varying} regions where much higher resolution is required
than elsewhere in the domain, particularly 
\revised{since the solution may contain
discontinuities or other localized features.}{when tracking a tsunami
propagating across the ocean}.
The \revised{Clawpack}{open source GeoClaw} software deals with 
this issue by using block-structured adaptive mesh 
refinement to selectively refine around propagating waves. For problems 
where only a target area of the total solution is of interest (e.g. one
coastal community), a method that 
allows identifying and refining the grid only in regions that influence this 
target area would significantly reduce the computational cost of finding a 
solution. 

In this work, we show that solving the time-dependent adjoint equation and using 
a suitable inner product with the forward solution allows more precise 
refinement of the relevant waves. 
\added{We present the adjoint methodology first in one space dimension for
illustration and in a broad context since it could also
be used in other adaptive software, and potentially for other tsunami
applications beyond adaptive refinement.}
\revised{
We present examples solving the shallow 
water equations in one and two dimensions. To perform these simulations,
the use of the adjoint method has been integrated into the adaptive 
mesh refinement strategy of the open source GeoClaw software.
We also present results that show}{
We then show how this
adjoint method has been integrated into the adaptive 
mesh refinement strategy of the open source GeoClaw
software and present tsunami modeling results showing}
that the accuracy of the solution is maintained 
and the computational time required is significantly reduced through the integration 
of the adjoint method into adaptive mesh refinement.

\end{abstract}


\section{Introduction}
\label{sec:intro}
Adjoint methods are often used in conjunction with the
numerical solution of differential equations for a variety of purposes,
such as computing sensitivities of the solution to input data, 
solving inverse problems, estimating errors in the solution,
or guiding the design of a computational grid to most efficiently compute
particular quantities of interest. 
We are exploring the use of adjoint methods in tsunami modeling by
incorporating them into the open source GeoClaw software that is widely
used for tsunami simulation, see e.g. 
\cite{BergerGeorgeLeVequeMandli2011,geoclaw,LeVequeGeorgeBerger2011}. 
We mention several other potential applications below, but
in this paper we focus primarily
on one particular use: guiding adaptive mesh
refinement (AMR) to efficiently capture the waves from a far-field tsunami that
will impact a particular ``target location,'' by which we mean
a specific location where we want to compare with available DART buoy or
tide gauge data, or a portion of the coastline where
we wish to compute inundation, for example.  

GeoClaw simulations often use 6 or 7 nested levels of refinement,
starting with a resolution of 1 or 2 degrees of latitude/longitude over
the entire computational domain.  This might be refined to 4-minute
or 1-minute resolution around the propagating waves, and then refined
to successively higher resolution around the target region, where
the finest grids may be 1/3 arcsecond for inundation studies.  
\revised{}{The target location where the highest levels of refinement is required
can be specified directly by the user.  In this paper we are concerned with
the question of how best to refine the ocean to capture the portions of the
propagating tsunami that will eventually impact this location, without
over-refining waves that will not.}

Refining the wave as it propagates across the ocean can be automated in
GeoClaw, by refining only in regions where the surface elevation differs
from sea level above some specified tolerance.  Every few time steps
recursive regridding is performed in which all such points on a given grid
are flagged for refinement to the next level and the flagged points are
clustered into rectangular patches using the algorithm of
\cite{BergerRigoutsos1991}.  Refinement ``regions'' can also be specified,
space-time subsets of the computational domain where
refinement above a certain level can be either required or forbidden.  This
is used to allow the finest levels of refinement only near the target
location. 
These AMR regions can also be used to induce the code to follow only
the waves of interest as the tsunami propagates across the ocean, but to
do so optimally often requires multiple attempts and careful examination of
how the solution is behaving, generally using coarser grid runs for guidance.
This manual guiding of AMR may also fail to capture some waves that are
important. For example, a portion of a tsunami wave may appear to be heading
away from the target location but might later reflect off a distant
shoreline or underwater feature, or edge waves may be excited that propagate
back and forth along
the continental shelf for hours after the primary wave has passed.

This challenge in tsunami modeling was the original motivation for 
our work on adjoint-based refinement, which we are also
incorporating into the more general Clawpack software (\cite{clawpack}),
which solves general hyperbolic partial differential equations
that arise in many wave propagation problems.  GeoClaw is based
on the AMRClaw branch of Clawpack, a more general code that implements
adaptive mesh refinement in both two and three space dimensions.
Other applications where adjoint-based refinement could be very
useful include earthquake simulation, for example, where the
desire might be to efficiently refine only the seismic waves that will reach
a particular location.

For a time-dependent partial differential equations such as the shallow
water equations, we generally wish to solve a ``forward problem'' in which 
initial data is specified (e.g. sea floor displacement due to an earthquake)
at some initial time $t_0$
and the problem is solved forward in time to find the effect at the target
location at some later time $t_f$.
The {\em adjoint equation} is a closely related partial differential
equation that must be solved {\em backwards} in time from the final time
$t_f$ to the initial time, as derived in \Sec{sec:modeling}. 
The data for the adjoint equation (specified at the final time) typically
approximates a delta function at the target location.  This spreads out into
waves as the adjoint equation is solved backward in time, with the same
\added{bathymetry-dependent} wave
speed $\sqrt{gh}$ as in the \revised{}{linearized} forward problem.  The key idea is
that at any intermediate time $t$ between $t_0$ and $t_f$, the only regions in space
where the forward solution could possibly reach the target location at time
$t_f$ are regions where the adjoint solution is nonzero.  Moreover by
computing a suitable inner product of the forward and adjoint solutions at
time $t$, it is possible to determine whether  the forward solution wave at
a given spatial point will actually reach the target location, or whether it
can be safely ignored.  This information can then be used to
decide whether or not to refine this spatial location in the forward solution.

In \Sec{sec:adjoint1} we briefly introduce the mathematical concept of an
adjoint equation.
Then in \Sec{sec:modeling} we derive the adjoint equation
for the linearized shallow water equations in one space dimension and
illustrate the main idea used in adjoint-based refinement in this simple
context before extending to the two-dimensional problem in
\Sec{sec:modeling_2d}. In \Sec{sec:geoclawex} we 
illustrate the use of this method and its efficiency on a tsunami modeling
problem using GeoClaw.

\section{Adjoint equation for a system of equations}\label{sec:adjoint1}
For readers not familiar with the concept of an adjoint equation, it may be
easiest to appreciate the power and limitations of this approach by first
considering the solution to an algebraic system of equations, beginning with
a linear system of the form $Ax = b$, where $A$ is an invertible $n\times n$
matrix, $b$ is a given vector with $n$ components, and the solution is
$x=A^{-1}b$.  In practice such a system is best solved by Gaussian
elimination, requiring ${\cal O}(n^3)$ operations.  Suppose that we are not
interested in the full solution $x$ but only in one component, say $x_k$. In
general we must still solve the full system to determine $x_k$. But now
suppose we want to do this for many different sets of data $b$, or that we
wish to determine the sensitivity of $x_k$ to changes in any component of
$b$.  In these situations the adjoint equation can be very useful since it
requires solving only a single system of equations rather than many systems.  

More generally, suppose that we only care about $J=\phi^Tx=\sum_{i=1}^n
\phi_i x_i$, where $\phi$ is
some specified vector with $n$ components.  In particular if $\phi$ is the
unit vector with $\phi_k=1$ and $\phi_i=0$ for $i\neq k$, then $\phi^Tx=x_k$,
the case considered in the previous paragraph.  The adjoint approach works
for more general $\phi$, i.e. when we only care about some scalar
quantity of interest
that can be defined as a linear functional applied to $x$.

For the linear system $Ax=b$, 
the adjoint equation is the linear system $A^T \hat x = \phi$, where the
vector $\phi$ is now used as the data on the right hand side and we solve
for $\hat x$, the adjoint solution.  The matrix $A^T$ is the transpose (also
called the adjoint) of the matrix $A$, with elements $(A^T)_{ij} = A_{ji}$.
This is also an $n\times n$ invertible
matrix so this problem has a unique solution $\hat x = A^{-T}\phi$.  
The matrix $A^{-T}$ is the inverse of the transpose, which agrees with the
transpose of the inverse.

The adjoint solution can now be used to compute $J(b)$, the value of the
functional $\phi^Tx$, where $x$ solves
$Ax=b$, by using elementary linear algebra:
\begin{equation}\label{adjlin1}
J(b) = \phi^Tx = \phi^T A^{-1}b = (A^{-T}\phi)^Tb = \hat x^T b.
\end{equation} 
Note that once we have solved the adjoint equation for $\hat x$, we can 
compute $J(b)$ for {\em any} data $b$ without solving additional linear
systems. We need only compute the inner product $\hat x^T b = \sum_{i=1}^n \hat
x_ib_i$, which requires only ${\cal O}(n)$ operations.  

Moreover, we can also compute the sensitivity of $J(b)$ to a change in any
component $b_i$ of the data.  Differentiating \cref{adjlin1} with respect to
$b_i$ shows that 
\begin{equation}\label{dJdbi}
\frac{\partial J(b)}{\partial b_i} = \hat x_i,
\end{equation} 
in other words the components of the adjoint solution are exactly the
sensitivities of $J$ to changes in the corresponding component of $b$.
We could have estimated the sensitivity ${\partial J(b)}/{\partial b_i}$
by varying $b_i$ slightly and solving a perturbed linear system, but we would
have had to solve $n$ such linear systems to estimate all the sensitivities.  The
adjoint equation computes them all simultaneously
through the solution of a single linear system.

In tsunami modeling we may wish to compute the sensitivity of the tsunami
observed at our target location to changes in the data, e.g. to changes in
the seafloor deformation if we are using a gradient-based optimization
algorithm to solve an inverse problem to match
observations (see \cite{Blaisea2013} for one such application to tsunami
source inversion).  Or we may want
to determine what potential source regions around the Pacific Rim give the
largest tsunami response at a particular target location (such as Pearl
Harbor, as considered in a study by \cite{Tang2006}).  Rather than
solving many forward problems, this can be determined with a single adjoint
solution. 

One limitation of the adjoint approach is that changing the target location
is analogous to changing the vector $\phi$ defining the quantity of interest
$J(b)$ in the linear system problem, and a new adjoint solution must be
computed for each location of interest.

Another limitation is that the adjoint approach is most easily applied to a
linear problem.  If we replace the linear system $Ax=b$ by a nonlinear
system of equations $f(x) = b$ that defines $x$ for data $b$, then we can
still use an adjoint approach to compute sensitivities of $J(b) = \phi^Tx$
to changes in $b$, but we must first 
linearize about a particular set of data $\bar b$
with solution $\bar x$ and can only compute sensitivities due to small
changes in $b$ around $\bar b$.  The adjoint equation then takes the form of a
linear system where the matrix $A$ is replaced by the Jacobian matrix of the
function $f$ evaluated at $\bar x$.

The GeoClaw software solves
the nonlinear shallow water equations, but in this paper we restrict our
attention to the above-mentioned application of tracking waves in the ocean
that will reach the target location.  Since a tsunami in the
ocean typically has an amplitude that is very small
compared to the ocean depth, these equations essentially reduce to the linear
shallow water equations and the adjoint equation linearized about the
ocean at rest is sufficient for our needs.  We will see in \Sec{sec:modeling} that
these adjoint equations take a very similar form to the linearized shallow water
equations, although with slightly different boundary conditions. 
If we wanted to compute sensitivies of the nonlinear 
onshore inundation to changes in data then we would have to
linearize about a particular forward solution.
In \Sec{sec:futurework} we make some additional
comments about extension to nonlinear problems.

Adjoint equations have been used computationally for many years in a variety 
of different fields, with wide ranging applications. A few examples 
include weather model tuning (\cite{Hall1986}), aerodynamics 
design optimization (e.g. \cite{GilesPierce2000,Jameson1988,KennedyMartins2013}), 
automobile aerodynamics (\cite{Othmer2004}), and geodynamics
(\cite{BungeHagelbergTravis2003}).  They have been used for
seismic inversion (e.g. \cite{AkcelikBirosGhattas2002,TrompTapeLie2005})
and recently also applied to tsunami inversion (\cite{Blaisea2013}).
The adjoint method 
has also been used for error estimation in the field of aerodynamics 
(\cite{BeckerRannacher2001}) and for general coupled time-dependent 
systems (\cite{AsnerTavenerKay2012}). 
Various solution methods have been combined with
adjoint approaches, including Monte Carlo (\cite{BuffoniCupini2001}), 
finite volume (\cite{Mishra2013}), finite element 
(\cite{AsnerTavenerKay2012}), and spectral-element 
(\cite{TrompTapeLie2005}) methods. 

It is also possible to use adjoint equations
to compute sensitivities of $J$ to changes in the input data.
This has led to the adjoint equations being utilized for system control 
in a wide variety of applications such as shallow-water wave control by
\cite{SandersKatopodes2000} and 
optimal control of free boundary problems by \cite{Marburger2012}.
This is also useful in solving inverse problems and potential applications of
this approach in tsunami modeling are being studied separately.

The adjoint method has also been used to guide adaptive mesh 
refinement, typically by estimating the error in the calculation 
and using that to determine how to adjust the grid.
Leveraging the adjoint problem to achieve this goal is not 
a new concept, and has been explored significantly for steady state problems 
where work has been done to guide AMR (e.g. 
\cite{PierceGiles2000, BeckerRannacher2001,VendittiDarmofal2000,
VendittiDarmofal2002,Park2004,VendittiDarmofal2003}),
and put error bounds on solution accuracy (e.g. \cite{GilesSuli2002}).
In the finite volume literature this approach is known as 
output-based mesh adaptation, although perhaps a clearer term would be 
\textit{adjoint-error} based mesh adaptation. 
Within the finite volume community output-based mesh 
refinement has begun to be used for unsteady
 problems. Specifically, temporal-only adaptation and space-time adaptation in the 
context of aerodynamics have been
explored by \cite{ManiMavriplis2007} and \cite{FlyntMavriplis2012} respectively, 
and work with the compressible Navier-Stokes equations has been done
for both static domains, by \cite{LuoFidkowski2011}, and deforming
domains, by \cite{KastFidkowski2013}.

\section{One-dimensional shallow water equations}\label{sec:modeling}
In one \revised{dimension these}{space dimension the shallow water equations}
take the form
\begin{subequations}\label{eq:swe_1d}
\begin{align}
h_t + ( hu)_x &= 0 \\
( hu )_t + ( hu^2 + \tfrac{1}{2}gh^2)_x &= -ghB_x.
\end{align}
\end{subequations}
Here, $u(x,t)$ is the depth-averaged velocity, 
$B(x)$ is the bottom surface elevation relative
to mean sea level, $g$ is the gravitational 
constant, and $h(x,t)$ is the fluid depth. We will use $\eta (x,t)$ to denote 
the water surface elevation,
\begin{align*}
\eta (x,t) = h(x,t) + B(x).
\end{align*}

The shallow water equations are a special case of a hyperbolic system
of equations,
\begin{equation}\label{eq:hypsys_1d}
q_t(x,t) + f(q)_x = \psi (q,x)
\end{equation} 
in one dimension and 
\begin{equation}\label{eq:hypsys}
q_t(x,y,t) + f(q)_x + g(q)_y= \psi (q,x,y)
\end{equation} 
in two dimensions, where $q$ is a vector of unknowns, $f(q)$ and $g(q)$ are the 
vectors of corresponding fluxes, and $\psi$ is a vector of source terms.
These appear in the study of numerous physical phenomena where wave motion is
important, and hence methods for numerically calculating solutions to
these systems of partial differential equations have broad applications over
multiple disciplines. 

As mentioned above, when tracking a tsunami in the ocean the nonlinear 
shallow water equations essentially reduce to the linear shallow water 
equations.  To linearize the shallow water equations about 
the ocean at rest, we begin by 
letting $\mu = hu$ represent the momentum and
noting that the momentum equation from \cref{eq:swe_1d} can be rewritten as
\begin{align*}
\mu_t + (hu^2)_x + gh(h+B)_x &= 0. \\
\end{align*}
Linearizing this equation as well as the continuity equation 
about a flat surface $\bar{\eta}$ and zero velocity 
$\bar{u} = 0$, with $\bar{h}(x) = \bar{\eta} - B(x)$ gives
\begin{subequations}\label{eq:swe_1d_linearized}
\begin{align}
\tilde{\eta}_t + \tilde{\mu}_x &= 0 \\
\tilde{\mu}_t + g \bar{h}(x)\tilde{\eta}_x &= 0 
\end{align}
\end{subequations}
for the perturbation $(\tilde{\eta}, \tilde{\mu})$ about 
$(\bar{\eta}, 0)$. 
Dropping tildes and setting 
\begin{align}\label{eq:SWE_Aq}
A(x) = \left[ \begin{matrix}
0 & 1 \\
g \bar{h}(x) & 0 
\end{matrix}\right],\hspace{0.1in}
q(x,t) = \left[\begin{matrix}
\eta \\ \mu
\end{matrix}\right],
\end{align}
gives us the system 
\begin{equation}\label{eq:hyp_linear}
q_t(x,t) + A(x)q_x(x,t) = 0.
\end{equation}

\revised{}{We now derive the adjoint equation for a linear hyperbolic
system of partial differential equations of the form \cref{eq:hyp_linear}
posed on an interval $a\leq x \leq b$ and over a time interval $t_0\leq t \leq
t_f$,}
subject to some known initial conditions, 
$q(x,t_0)$, and some boundary conditions at $x=a$ and $x=b$. 
In the linearized shallow water case $A(x)$ and $q(x,t)$ 
are given by \cref{eq:SWE_Aq},
although the analysis shown below applies more generally to any 
time dependent hyperbolic system of equations.
In \Sec{sec:adjoint1} we considered the case where we care about 
$J = \phi^Tx$ where $x$ was the solution to an algebraic
system of equations, now suppose
we are interested in calculating the value of a functional
\begin{equation}\label{eq:J_general}
J = \int_a^b \varphi^T (x) q(x,t_f) dx
\end{equation} 
for some given $\varphi (x)$.
For example, if $\varphi (x) = \delta (x - x_0)$ then $J = q(x_0,t_f)$ is the
solution value at the point $x = x_0$ at the final time $t_f$. 
This is the situation we consider in this paper, with the delta function
smeared out around the region of interest for the computational approach.

If $\hat{q}(x,t)$ is any other appropriately sized 
vector of functions then multiplying this by
\cref{eq:hyp_linear} and integrating in both space and time yields
\begin{equation}\label{eq:q_integral}
\int_a^b \int_{t_0}^{t_f}
\hat{q}^T(x,t)\left(q_{t}(x,t)+A(x)q_{x}(x,t)\right)dx\,dt = 0
\end{equation}
for any time $t_0 < t_f$.
Then integrating by parts \revised{twice}{in space and then in time}
yields the equation 
\begin{equation}\label{eq:intbyparts}
\int_a^b   \left.\hat{q}^Tq\right|^{t_f}_{t_0}dx 
+ \int_{t_0}^{t_f} \left.\hat{q}^TAq\right|^{b}_{a}dt 
- \int_{t_0}^{t_f} \int_a^b  q^T\left(\hat{q}_{t} +
\left(A^T\hat{q}\right)_{x}\right)dx\,dt = 0.
\end{equation} 
By defining the adjoint equation, 
\begin{equation}\label{adjoint1}
\hat{q}_{t}(x,t) + (A^T(x)\hat{q}(x,t))_{x} = 0,
\end{equation}
setting $\hat{q}(x,t_f) = \varphi (x)$,
and selecting the appropriate boundary conditions for $\hat{q}(x,t)$ 
such that the integral in time vanishes (which varies based on the 
specific system being considered), we can eliminate all
terms from \cref{eq:intbyparts} except the first term, to obtain 
\begin{equation}\label{eq:q_equality}
\int_a^b \hat{q}^T(x,t_f)q(x,t_f) dx = \int_a^b
\hat{q}^T(x,t_0)q(x,t_0)dx.
\end{equation} 
Therefore, the integral of the inner product
between $\hat{q}$ and $q$ at the final time is equal 
to the integral at the initial time $t_0$:
\begin{equation}
J = \int_a^b\hat{q}^T(x,t_0)q(x,t_0)dx.\label{eq:J}
\end{equation}
Note that we can replace $t_0$ in \cref{eq:q_integral} with any 
$t$ so long as $t_0 \leq t \leq t_f$, which would yield \cref{eq:J} 
with $t_0$ replaced by $t$. 
 From this we observe that the locations where the 
 magnitude of the inner product
$\hat{q}(x,t)^Tq(x,t)$ is large, for any $t$ with
$t_0 \leq t \leq t_f$, are the areas that will have a
significant effect on the inner product $J$.
These are the \revised{areas where the solution should be refined}{candidate
areas for refinement} at time $t$. 
To make use of this, we must first solve the adjoint equation
\cref{adjoint1} for $\hat q(x,t)$.
Note, however, that this requires using 
``initial'' data $\hat{q}(x,t_f)$, so
the adjoint problem must be solved backward in
time.
\revised{}{In \Sec{sec:amradj} we will discuss in more detail the manner in
which this is done.  }

First, we present a one-dimensional example that illustrates how the
waves from the forward and adjoint equations propagate and can be combined
to identify the waves that will reach a point of interest.
Viewing this first in one dimension has the advantage that we can 
easily view the full time-history of
waves in the forward and adjoint equations together in single plots in the
$x$-$t$ plane.  The figures presented in this section were computed using
Clawpack on a very fine grid to generate illustrations, but the particular
numerical method is immaterial here.  We are not testing the adjoint approach
to AMR since we do not perform any adaptive refinement in this case.  In 
\Sec{sec:modeling_2d} below we will present numerical results using GeoClaw with
adaptive refinement and further discuss the numerical methods used there.

As an example we use the one-dimensional linearized shallow water equations
and its adjoint on the domain $0 \leq x \leq 400$ km. 
We choose simplified topography that is piecewise constant,
with depth 4000 m for $50 \leq x \leq 400$ and shallower depth 200 m for
$0 \leq x \leq 50$.  The step discontinuity from deep ocean to ``continental
shelf'' is chosen so that waves remain localized and the plots are easy to
interpret.  For the same reason we use reflecting boundary conditions at each
boundary rather than a more realistic shore: 
\begin{align*}
u(0,t) = 0, \hspace{0.1in} u(400,t) &= 0 &&t  \geq t_0.
\end{align*}

\begin{figure}[h!]
\begin{minipage}[b]{0.5\linewidth}
\includegraphics[width=\textwidth]{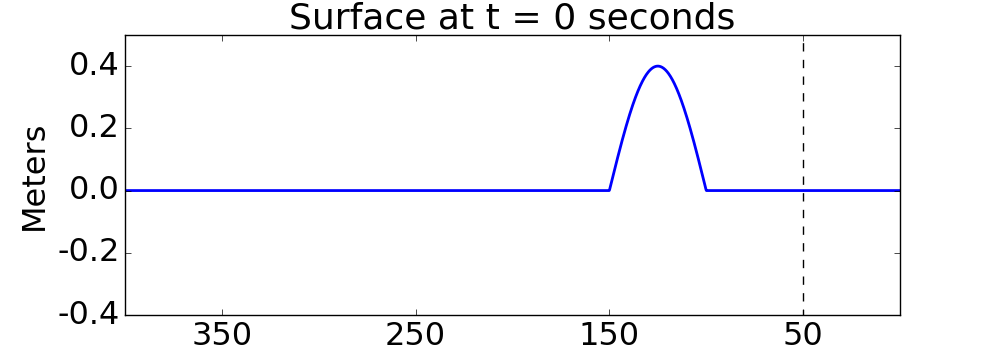}
\end{minipage}
\vspace{0.2cm}
\begin{minipage}[b]{0.5\linewidth}
\includegraphics[width=\textwidth]{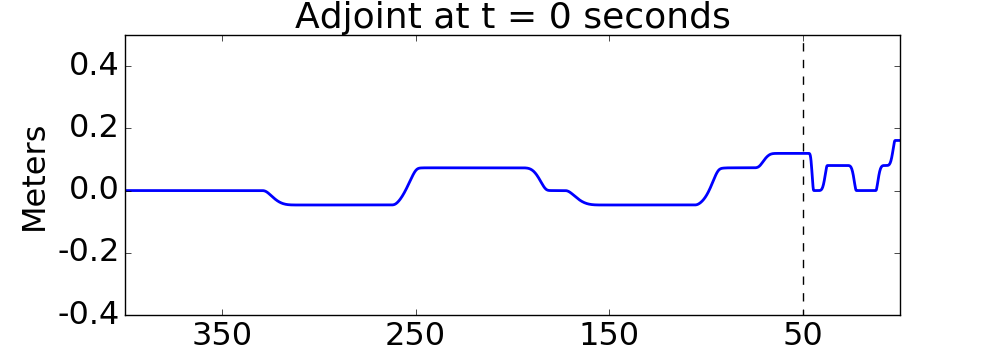}
\end{minipage}
\begin{minipage}[b]{0.5\linewidth}
\includegraphics[width=\textwidth]{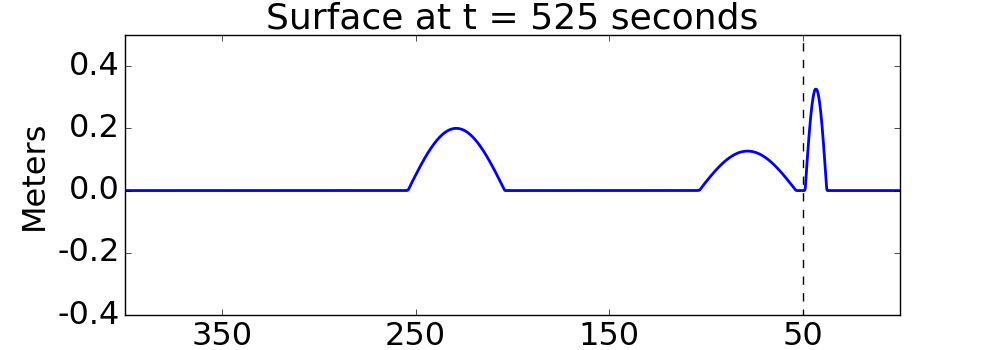}
\end{minipage}
\vspace{0.2cm}
\begin{minipage}[b]{0.5\linewidth}
\includegraphics[width=\textwidth]{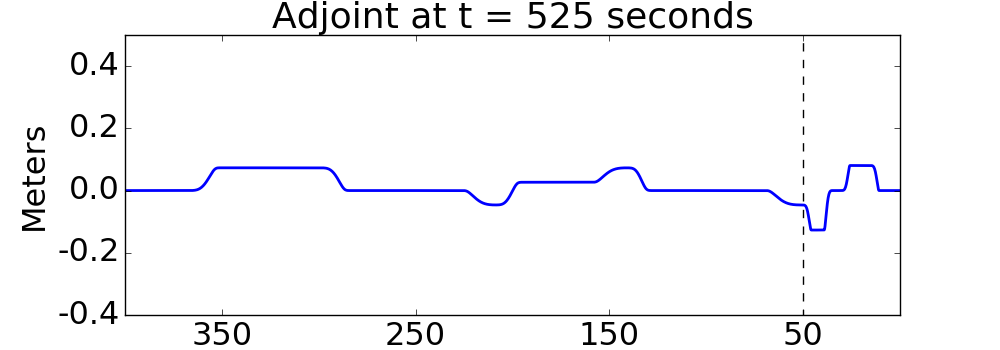}
\end{minipage}
\begin{minipage}[b]{0.5\linewidth}
\includegraphics[width=\textwidth]{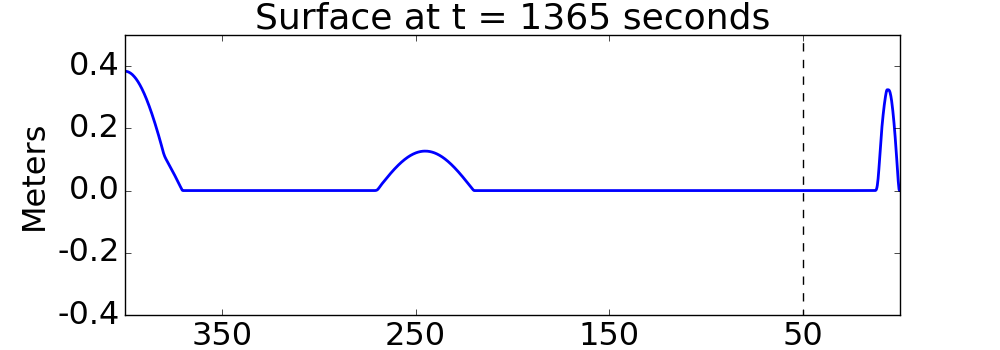}
\end{minipage}
\vspace{0.2cm}
\begin{minipage}[b]{0.5\linewidth}
\includegraphics[width=\textwidth]{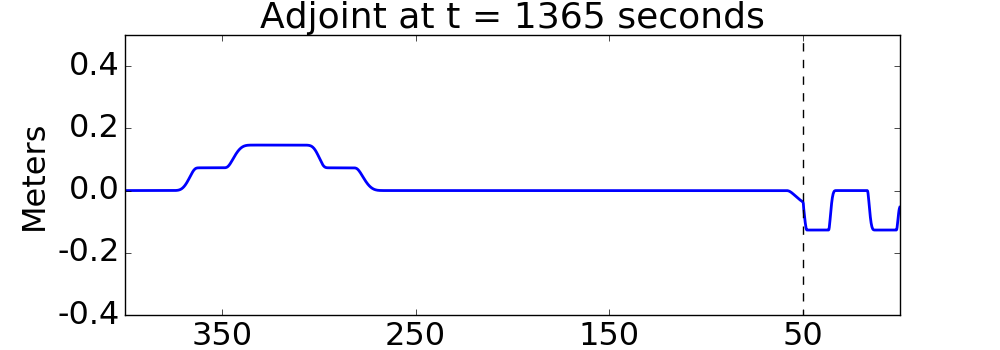}
\end{minipage}
\begin{minipage}[b]{0.5\linewidth}
\includegraphics[width=\textwidth]{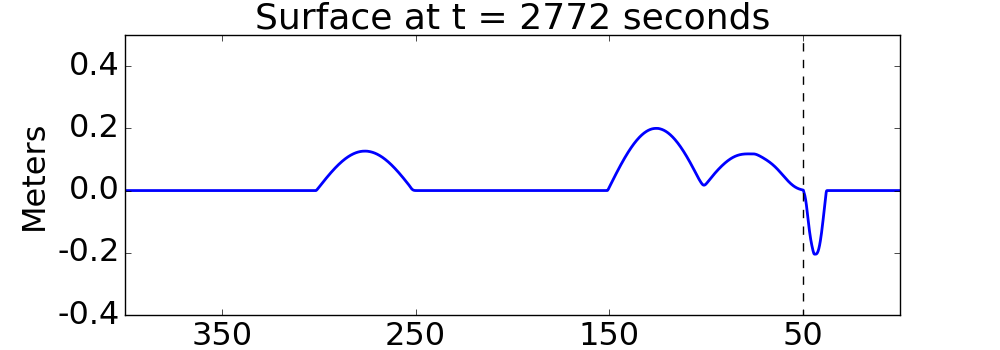}
\end{minipage}
\vspace{0.2cm}
\begin{minipage}[b]{0.5\linewidth}
\includegraphics[width=\textwidth]{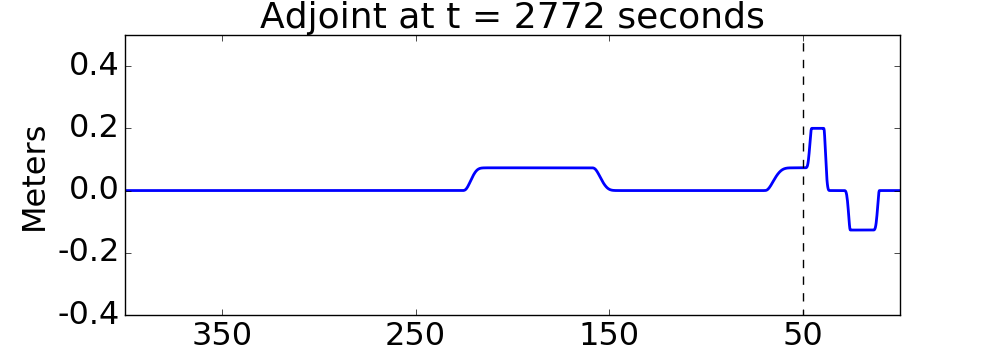}
\end{minipage}
\begin{minipage}[b]{0.5\linewidth}
\includegraphics[width=\textwidth]{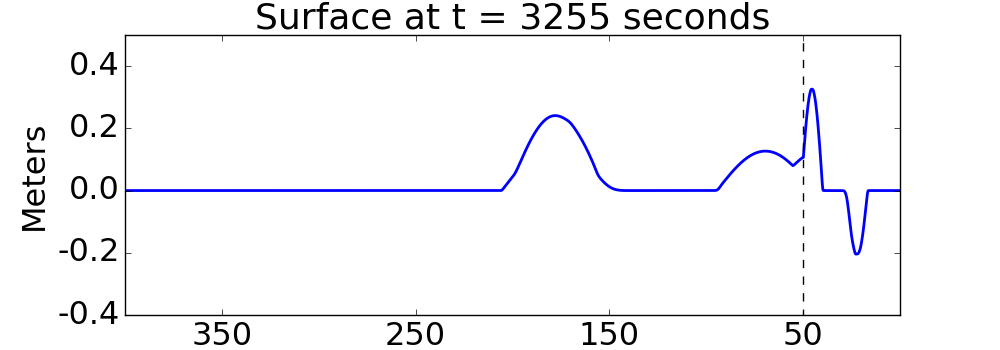}
\end{minipage}
\vspace{0.2cm}
\begin{minipage}[b]{0.5\linewidth}
\includegraphics[width=\textwidth]{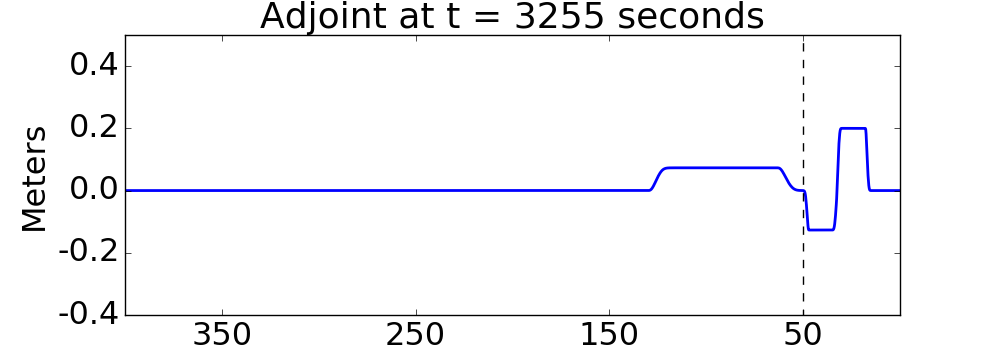}
\end{minipage}
\begin{minipage}[b]{0.5\linewidth}
\includegraphics[width=\textwidth]{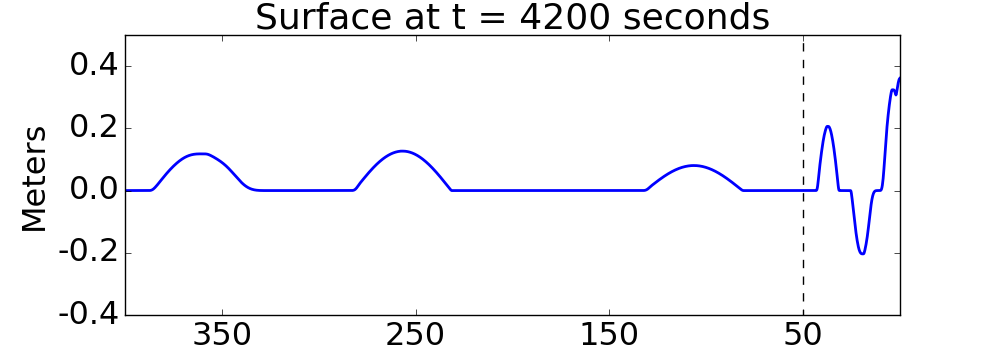}
\end{minipage}
\begin{minipage}[b]{0.5\linewidth}
\includegraphics[width=\textwidth]{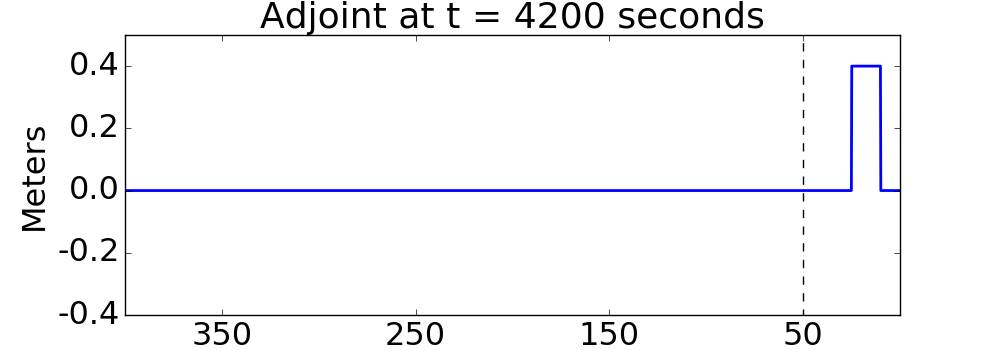}
\end{minipage}
 \caption{%
  An idealized tsunami interacting with a step discontinuity in bathymetry, 
  for both the forward and the adjoint problems. The 
  dashed line indicated the location of the discontinuity, $50$ km offshore.}
\label{fig:1dSWE}
\end{figure}

As initial data for $q(x,t)$ we introduce a hump of stationary water 
with a height of $0.4$ m centered $125$ km off shore. 
\Fig{fig:1dSWE} shows the resulting waves interacting with 
the boundaries and the discontinuity in bathymetry, at
the location indicated by the dashed line. 
As time progresses the hump splits into 
equal \revised{left}{right}-going and \revised{right}{left}-going waves heading out towards the ocean 
and towards the continental shelf, respectively. When the \revised{right}{left}-going 
wave encounters the continental shelf it splits into a reflected and a 
transmitted wave. When any waves, e.g. either the original 
\revised{left}{right}-going wave or the newly transmitted wave, encounter the wall on either side 
of the domain the waves are reflected back in the other direction. This 
interplay between reflection and transmission of the waves at the 
bathymetry discontinuity and the reflection of waves at the two boundaries
leads to a \revised{more complex wave pattern than would be observed if only 
the initial wave were considered}{complex wave pattern}.

For the adjoint solution, suppose that 
we are interested in the accurate estimation of the 
\revised{pressure}{surface elevation} in the 
interval between $10$ and $25$ km offshore, perhaps because that 
is the location of \revised{a gauge}{gauges} with which we wish to compare our results.

\subsection{Single Point In Time}\label{sec:1d_timepoint}
\revised{}{Initially suppose we are only interested in one particular time,
say $t_f=4200$ seconds.} 
Setting $J = \int_{10}^{25}\eta(x,t_f)\,dx$,
the problem then requires that
\begin{align}
\varphi (x) = \left[ \begin{matrix}
I (x) \\ 0
\end{matrix}
\right], \label{eq:phi_1d}
\end{align}
where
\begin{align}
I (x) = \left\{
     \begin{array}{lr}
       1 & \hspace{0.3in}\textnormal{if } 10 < x < 25\\
       0 & \textnormal{otherwise.}\hspace{0.21in}
     \end{array}
   \right. \label{eq:delta_1d}
\end{align}
Define
\begin{align*}
\hat{q}(x,t_f) = \left[ \begin{matrix}
\hat{\eta}(x,t_f) \\ \hat{u}(x,t_f)
\end{matrix}\right] = \varphi(x),
\end{align*}
and note that \cref{eq:q_integral} holds for this problem.
If we define the adjoint problem by
\begin{align*}
\hat{q}_{t} + \left(A^T(x)\hat{q}\right)_{x} &= 0 &&x \in [0,400],
\hspace{0.1in}t_0 \leq t \leq t_f\\
\hat{u}(0,t) = 0, \hspace{0.1in} \hat{u}(400,t) &= 0 && t_0 \leq t \leq t_f,
\end{align*}
then \cref{eq:q_equality} holds,
which is the expression that allows us to use the inner product of the
adjoint and forward problems at each time step to determine what regions will
influence the point of interest at the final time. 

As the ``initial'' data for $\hat{q}(x,t_f) = \varphi (x)$ 
we have a square pulse in water height, which was described above in
Equations \cref{eq:phi_1d} and \cref{eq:delta_1d} at the final time.
As time progresses backwards, the pulse splits into equal left-going and right-going
waves which interact with the walls and the bathymetry discontinuity giving both reflected
and transmitted waves. 
\Fig{fig:1dSWE} also shows the results of solving this adjoint
problem.

\begin{figure}[h!]
 \centering
  \includegraphics[width=0.95\textwidth]{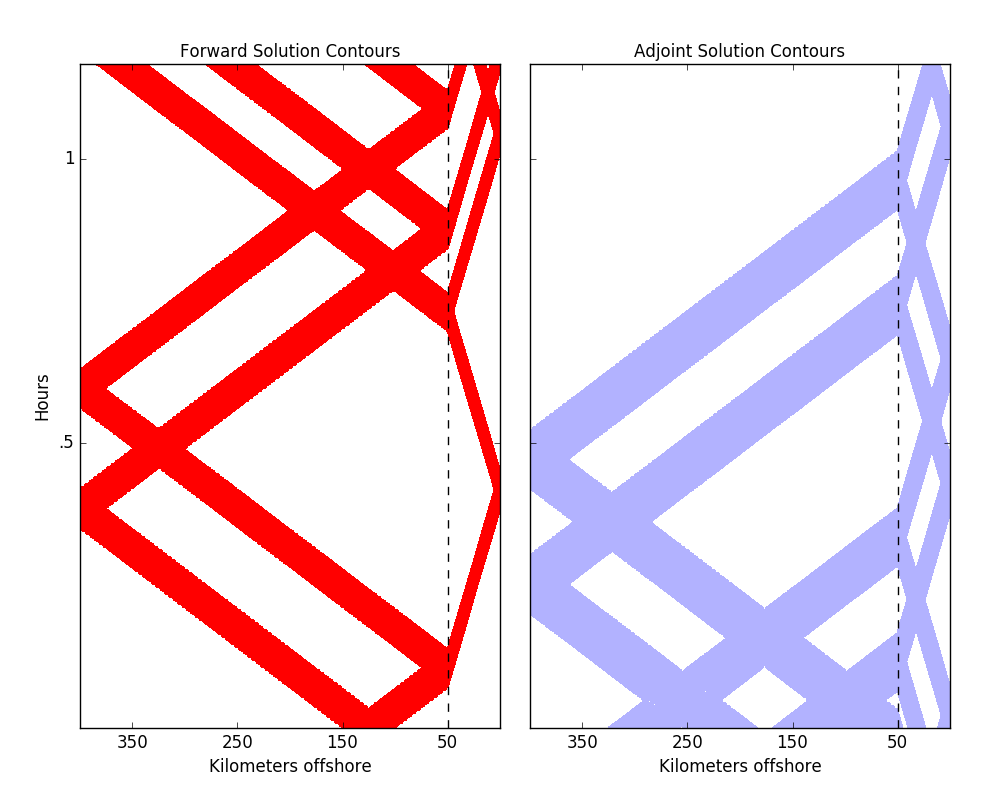}
 \caption{%
  \revised{Contour plot of the magnitude of the waves relative to mean sea
level.}{}
  On the left, the locations where the magnitude of the surface elevation
 $\eta$ in the forward problem is greater than or
 equal to $0.1$ meters above mean sea level are shown in red. 
 On the right, the locations where the magnitude of $\hat\eta$ in the
 adjoint solution is greater than or equal to $0.1$ are shown in blue. 
 \revised{}{The time axis is the same for both plots.}}
\label{fig:swe1d_xtplane}
\end{figure}

To better visualize how the waves are moving through the domain, it is helpful
to look at the data in the $x$-$t$ plane as shown in \Fig{fig:swe1d_xtplane}. 
For \Fig{fig:swe1d_xtplane}, the
horizontal axis is the position, $x$, and the vertical axis is time. The left
plot shows in red the locations where the magnitude of height of 
the tsunami in the forward problem is greater than or
equal to $0.1$ meters above mean sea level. 
On the right side of \Fig{fig:swe1d_xtplane} we show the adjoint
solution that is computed starting with the square wave data at the
final time $t_f$, indicating in blue the regions where the first
component $\hat \eta$ is greater than or equal to $0.1$.

\begin{figure}[h!]
  \centering
    \includegraphics[width=0.95\textwidth]{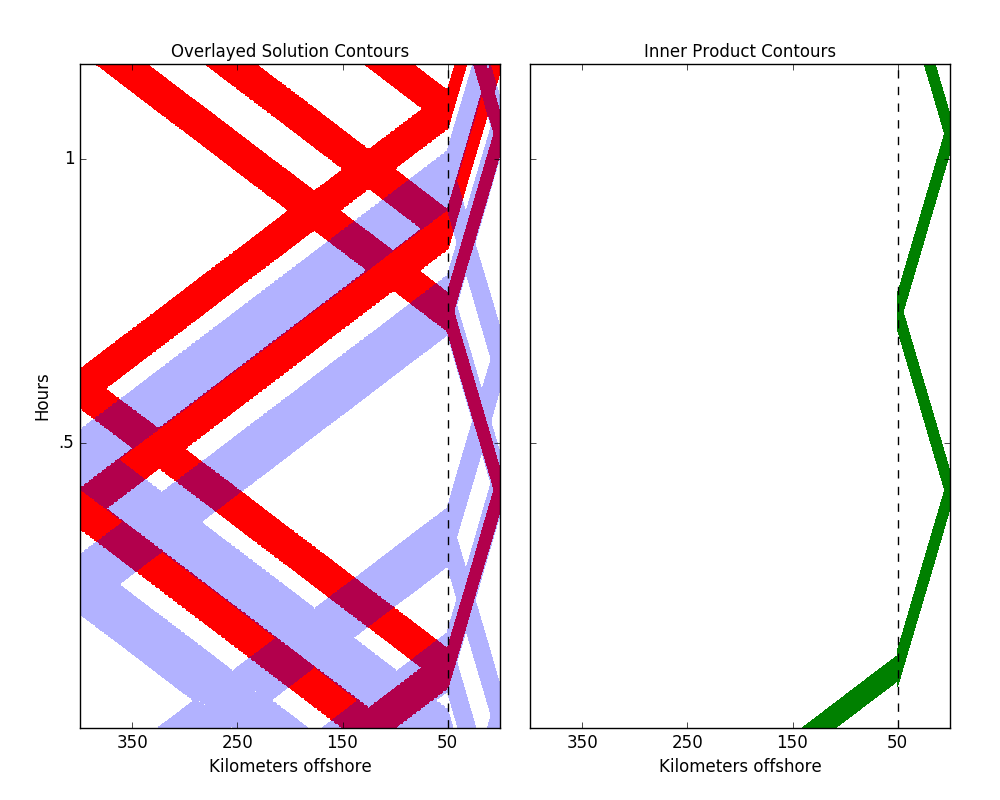}
    \caption{The \revised{contour}{} plots of the 
    magnitudes of $\eta$ and $\hat\eta$ from \Fig{fig:swe1d_xtplane}
    are overlayed in the left figure. The \revised{contour plot of
    the}{region where the} inner
    product of the two vectors $q(x,t)$ and $\hat{q}(x,t)$ \added{is above a
    threshold} is shown on the right, 
    \revised{indicating where the forward problem should be
    refined.}{picking out only the wave that reaches the target 
    location at $10 \leq x \leq 25$ the final time.}
    The time axis is the same for both plots.}
    \label{fig:swe1d_timepoint}
\end{figure}

On the left
side of \Fig{fig:swe1d_timepoint} we have overlayed the
magnitudes of the water heights for the forward and 
adjoint solutions, and 
viewing the data in the $x$-$t$ plane makes it fairly clear which parts 
of the wave from the forward
solution actually effect our region of interest at the final time. 

\Fig{fig:swe1d_timepoint}
also shows, in green on the right, the locations where the 
magnitude of the inner product between the
forward and adjoint solution is greater than or equal to 0.1 as time
progresses. At each time step this is clearly identifying the regions in the
computational domain that will affect our region of 
interest at the final time. If we were using adaptive mesh
refinement, these areas, identified by where the 
magnitude of $\hat{q}^T(x,t)q(x,t)$ 
exceeds some tolerance, are 
\revised{precisely the areas we would flag for}{the areas we should consider
for adaptive} 
refinement. Note that a mesh refinement strategy 
based \revised{on wherever}{only on the areas where} 
the magnitude of $\eta (x,t)$ is large would result in refinement of
many areas in the computational domain that will have no effect on our area of
interest at the final time (all the red regions in the left plot of 
\Fig{fig:swe1d_timepoint}).

\added{It is important to note that in the left side of
\Fig{fig:swe1d_timepoint} there are places where red and blue waves overlap
that do not show up in the plot on the right side.  These are areas where a
wave moving in one direction in the forward solution crosses a wave moving
in the other direction in the adjoint solution.  Even though both $q(x,t)$
and $\hat q(x,t)$ are nonzero vectors in these regions, the inner product of
the two is equal to zero.  This is easily verified by computing the
eigenvectors of the coefficient matrices $A$ and $A^T$ defining the
hyperbolic problems for $q$ and $\hat q$ (see, for example
\cite{LeVeque2002}).   The eigenvectors of
$A$ are $[1,~\pm c]^T$, where $c = \sqrt{g\bar h}$ is the wave speed.  
Hence the solution $q$ in a purely left-going wave is proportional to the
vector $[1,~-c]^T$, while in a purely right-going wave is proportional to
$[1,~+c]^T$.  The eigenvectors of $A^T$, on the other hand, are $[\pm
c,~1]^T$  and each of these is orthogonal to one of the eigenvectors of $A$
(and hence the inner product is zero for crossing waves).  It is only when
the waves are aligned on the left side of \Fig{fig:swe1d_timepoint} that the
inner product is seen to be nonzero on the right side of the figure. 
This further illustrates the power of the adjoint approach to identify only
the waves that will reach the target location.}

\subsection{Time Range}\label{sec:1d_timerange}

\begin{figure}[h!]
  \centering
    \includegraphics[width=0.95\textwidth]{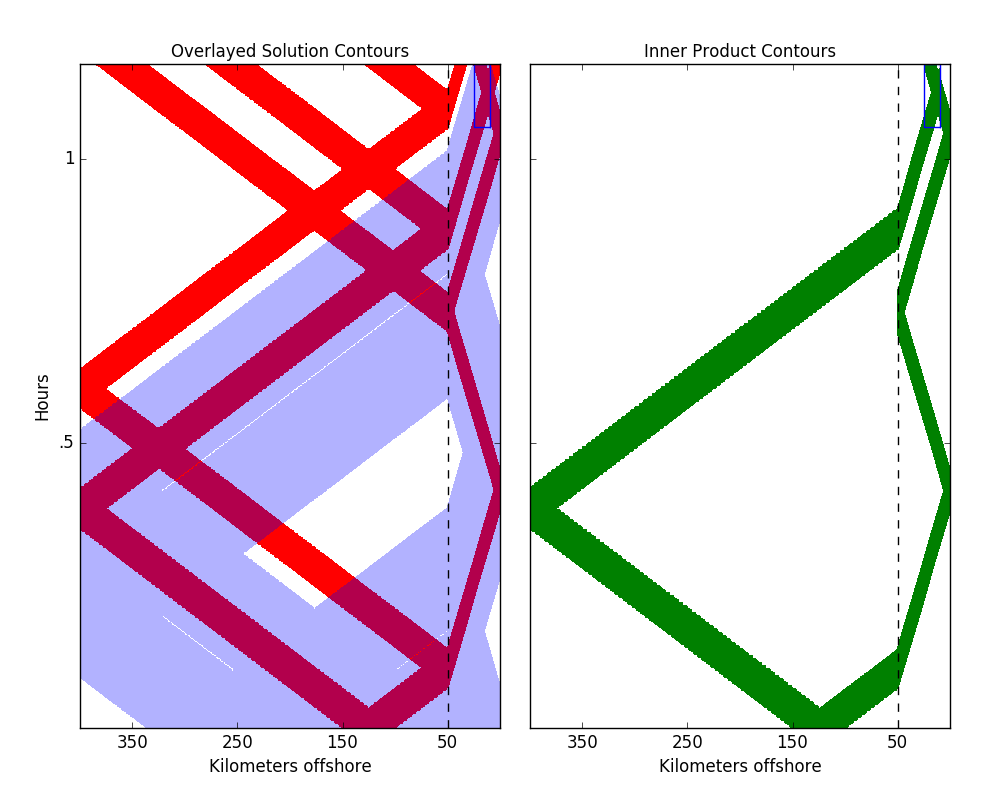}
    \caption{Computed results for one-dimensional shallow water equations in the case
    when the solution $q(x,t)$ is desired at the 
    \added{target location $10 \leq x \leq 25$ km
    over a time interval from 3800 to 4200 seconds,}
    as indicated by the box in each $x$--$t$ plot.
    Left: \revised{Contour plot for }{Plots showing}
    $\eta (x,t)$ and the adjoint solutions $\hat{\eta }(x,\tau)$ for 
    shifted values of $\tau$ as described in the text.
    Right: Regions where the maximum inner product over the given 
    time range exceeds the threshold, 
    showing the additional wave that \revised{must be refined.}{}reaches the
target region relative to the case shown in \Fig{fig:swe1d_timepoint}.
    The time axis is the same for both plots.}
    \label{fig:1dAcoustics_timerange}
\end{figure}

\added{
In practice we are rarely interested in the tsunami amplitude at some
location at only a single time, we are more likely to be interested in the
solution over some time range, often over the entire simulation time.
But note that once
we have computed the adjoint solution going backward from time $t_f$,
we immediately know the adjoint solution starting at some earlier time $\bar
t < t_f$: it is simply the same solution translated earlier in time by $t_f
- \bar t$, since the linearized adjoint equation is autonomous in
time. The formulas will be made more precise in the next section, but in
terms of \Fig{fig:swe1d_timepoint} we can think of moving the data for the
adjoint to an earlier time as simply translating the blue
solution in the left figure downward by this time increment.  
\Fig{fig:1dAcoustics_timerange} illustrates the result if we consider
setting $\bar t$ to all possible values between 
$t_s = 3800$ seconds and $t_f = 4200$ seconds.  The small blue rectangle in
the top corner of this plot shows the region in space-time that we are now  
concentrating on as our target space-time region, 
and the green portions in the figure on the right of 
\Fig{fig:1dAcoustics_timerange} show where the inner product of the forward
solution and any of these translated adjoint solutions is above the
threshold of 0.1.  We clearly see that, relative to
\Fig{fig:swe1d_timepoint}, there is an additional wave that must be tracked
since it arrives at the target spatial region during the target time period.

If we were interested in any wave that can arrive in the target region over
the full simulation time $0 \leq t \leq t_f$, then we would extend this by
looking at all possible downward translates of the adjoint solution.  
In this particular example, this would not identify any additional waves
beyond those already found, since the other waves would not arrive at the
target location until times beyond $t_f$.
}

\section{Combining mesh refinement and the adjoint problem}\label{sec:amradj}

\revised{
By developing a strategy for taking advantage of the information
provided by adjoint methods and incorporating it into an AMR algorithm, it is
possible to significantly reduce the computational time required to find the
solution. 
}{We now discuss in more detail how the adjoint solution might be used in
guiding adaptive mesh refinement, still focusing on the one-dimensional case
presented above.}

Before solving the forward problem, we would first solve the adjoint
equation \cref{adjoint1} backward in time. 
\revised{Clawpack}{Since GeoClaw and other numerical software} is designed to
solve equations forward in time, we consider the function
\begin{align*}
\tilde{q}(x,t) &= \hat{q}(x,t_f - t).
\end{align*}
This gives us the new problem
\begin{align*}
&\tilde{q}_t - \left( A^T(x) \tilde{q}\right)_x = 0 
&&x \in [a,b], \hspace{0.1in}t > 0\\
&\tilde{q}(a,t) = \hat{q}(a,t_f - t)  &&0 \leq t \leq t_f - t_0 \\
&\tilde{q}(b,t) = \hat{q}(b,t_f - t) &&0 \leq t \leq t_f - t_0
\end{align*}
with initial condition $\tilde{q}(x,0) = \varphi (x) $. This
problem is then solved \revised{using the Clawpack software.}{forward in time.} Snapshots of this solution
are saved at regular time intervals, $t_0, t_1, \cdots, t_N$, \revised{After the
solution is calculated,}{from which}
snapshots of the adjoint solution are retrieved by
simply setting
\begin{align*}
\hat{q}(x,t - t_n) = \tilde{q}(x,t_n)
\end{align*}
for $n = 0, 1, \cdots, N$. 

With the adjoint solution in hand, we now turn to the forward problem. 
\revised{As refinement occurs in space for the forward problem, 
maintaining the stability
of the finite volume method requires that refinement must also occur in time.
Therefore, as the forward solution is refined, solution data for the adjoint
problem is no longer available at the corresponding times (since only
snapshots at regular time intervals of the adjoint solutions were saved) or
locations (since the adjoint solution was calculated on a coarse grid).}{
In our approach to GeoClaw, described more fully in the two-dimensional case
in \Sec{sec:modeling_2d}, we solve the adjoint equation on a fixed spatial resolution
and output at discrete times.  To flag cells for refinement as we now solve
the forward problem, we will generally need to estimate the adjoint solution
on finer grids in both space and time.}
To address this issue, the solution for the adjoint problem at the necessary
\revised{locations is approximated using bilinear interpolation from the data present
on the coarser grid}{spatial locations is approximated using linear (or bilinear
in 2D) interpolation from the data present
on the coarser grid defined by the snapshots}.

\subsection{Single Point In Time}\label{sec:timepoint_theory}
If we are interested in the solution of the forward problem at our target area
\added{only} at the final time, then
when solving the forward problem we take the inner product
between the current time step in the forward problem and the 
time step in the adjoint problem in order to determine which areas in the
forward wave are going to impact the region of interest.
\revised{To be conservative, when}{When} considering the forward problem
at time $t$ \revised{with}{that lies between two snapshot times of the 
adjoint, say}
\begin{align*}
t_n \leq t \leq t_{n+1},
\end{align*}
 both $\hat{q}(x,t_n)$ and $\hat{q}(x,t_{n+1})$ are taken
into account.
\added{Rather than interpolating in time, }
since we wish to refine any parts of the domain that will have a significant impact 
on the final solution, we \added{take a more conservative approach and} 
refine wherever the magnitude of the inner product 
\begin{equation}\label{eq:ip_timepoint}
\max\limits_{\tau = t_n, t_{n+1}} 
 \left|\hat q^T(x, \tau)q(x,t)\right|
\end{equation}
is above some tolerance.
\added{Note that we must save snapshots of the adjoint at sufficiently dense
output times for this to be sufficient.}
 These areas are then
flagged for refinement, and the next time step is taken.

\subsection{Time Range}\label{sec:timerange_theory}
Suppose that we are instead interested in the accurate estimation 
of the forward problem in
same location for some time range $t_s \leq t \leq t_f$, where 
$t_0 \leq t_s \leq t_f$.
Define $\hat q(x,t; \overline{t})$ as the adjoint based
on data $\hat q(x,\overline{t}) = \varphi(x)$.  Then for each $\hat t$ in the interval
$[t_s,t_f]$, we need to consider the inner product of $q(x,t)$ with $\hat q(x,t;
\hat t)$. Note that since the adjoint is autonomous 
in time, $\hat q(x,t; \hat t) = \hat q(x, t_f-\hat t+t; t_f)$. Therefore, we must 
consider the inner product 
\begin{equation*}
\hat q^T(x, t_f-\hat t+t; t_f)q(x,t)
\end{equation*}
for $\hat t \in [t_s,t_f]$. Since we are in fact only concerned when 
the magnitude of this 
inner product is greater than some tolerance, we can simply consider
\begin{equation*}
\max\limits_{t_s \leq \hat t \leq t_f} 
 \left|\hat q^T(x, t_f-\hat t+t; t_f)q(x,t)\right|
\end{equation*}
and refine when this maximum is above the given tolerance. 
Define $\tau = t_f-\hat t+t$. Then this maximum can be 
rewritten as
\begin{equation}\label{eq:ip_timerange}
\max\limits_{T\leq \tau \leq t} 
\left|\hat{q}^T(x,\tau;t_f)q(x,t)\right|
\end{equation}
where $T = \min(t+t_f-t_s,~ t_f)$.

We now drop the cumbersum notation $\hat q(x,t; t_f)$ in favor of the simpler 
$\hat{q}(x,t)$ with the understanding that the adjoint is based on the data
$\hat q(x,t_f) = \varphi(x)$.  Note that we still only need to solve one
adjoint problem in this case, we simply use it over a larger time range in
evaluating \cref{eq:ip_timerange} than in \cref{eq:ip_timepoint}.

\section{Two-dimensional shallow water equations }\label{sec:modeling_2d}

In two space dimension the shallow water equations take the form
\begin{subequations}\label{eq:swe}
\begin{align}
h_t + ( hu)_x + ( hv)_y &= 0 \\
( hu )_t + ( hu^2 + \tfrac{1}{2}gh^2)_x + 
( huv)_y &= -ghB_x \\
( hv)_t + ( huv)_x + ( hv^2 + \tfrac{1}{2}gh^2)_y 
&= -ghB_y.
\end{align}
\end{subequations}
Here, $u(x,y,t)$ and $v(x,y,t)$ are the depth-averaged velocities in the two 
horizontal directions and $B(x,y)$ is the bottom surface elevation relative
to mean sea level. Now the water surface elevation is given by
\begin{align*}
\eta (x,y,t) = h(x,y,t) + B(x,y).
\end{align*} 
Linearizing \revised{the two dimensional shallow water}{these} equations 
about an ocean at rest, similar to what was done in \Sec{sec:modeling} for 
the one dimensional case, gives
\begin{align*}
\tilde{\eta}_t + \tilde{\mu}_x + \tilde{\gamma}_y &= 0 \\
\tilde{\mu}_t + g \bar{h}(x,y)\tilde{\eta}_x &= 0 \\
\tilde{\gamma}_t + g\bar{h}(x,y)\tilde{\eta}_y &= 0 
\end{align*}
for the perturbation $(\tilde{\eta}, \tilde{\mu}, \tilde{\gamma})$ about 
$(\bar{\eta}, 0, 0)$. 
Dropping tildes and setting 
\begin{align*}
A_1(x,y) = \left[ \begin{matrix}
0 & 1 & 0 \\
g \bar{h}(x,y) & 0 & 0 \\
0 & 0 & 0
\end{matrix}\right],\hspace{0.1in}
A_2(x,y) = \left[ \begin{matrix}
0 & 0 & 1 \\
0 & 0 & 0 \\
g \bar{h}(x,y) & 0 & 0
\end{matrix}\right], \hspace{0.1in}
q(x,y,t) = \left[\begin{matrix}
\eta \\ \mu \\ \gamma
\end{matrix}\right],
\end{align*}
gives us the system 
\begin{equation}\label{eq:2d_hyperbolicsystem}
q_t(x,y,t) + A_1(x,y)q_x(x,y,t) + A_2(x,y)q_y(x,y,t) = 0. 
\end{equation}

\revised{}{For a general system of the form \cref{eq:2d_hyperbolicsystem} we now derive the adjoint equation posed on an interval $a \leq x \leq b$, $\alpha \leq y \leq \beta$ and over a time interval $t_0 \leq t \leq t_f$, subject to some known initial conditions $q(x,y,t_0)$ and some boundary conditions at $x = a$, $x = b$, $y = \alpha$ and $y = \beta$. If $\hat{q}(x,y,t)$ is an appropriately sized vector of functions then}
note that
\begin{align*}
\int_{t_0}^{t_f}\int_a^b\int_{\alpha}^{\beta}  \hat{q}^T\left( q_t + A_1(x,y)q_x +
A_2(x,y)q_y\right) dy\,dx\,dt = 0.
\end{align*}
Following the same basic steps we used in one dimension to go from 
\cref{eq:q_integral} to \cref{eq:intbyparts}, integrating by parts 
yields the equation
\begin{align}\label{eq:intbyparts2d}
\int_a^b  \int_{\alpha}^{\beta}  \hat{q}^Tq|^{t_f}_{t_0}&dy\,dx +
\int_{t_0}^{t_f}\int_{\alpha}^{\beta}  \hat{q}^TA_1(x,y)q|^{b}_{a}dy\,dt 
+ \int_{t_0}^{t_f}\int_{a}^{b}  \hat{q}^TA_2(x,y)q|^{\beta}_{\alpha}dx\,dt \nonumber \\
&- \int_{t_0}^{t_f}\int_a^b\int_{\alpha}^{\beta}  q^T\left(\hat{q}_{t} +
\left(A_1^T(x,y)\hat{q}\right)_{x} + \left(A_2^T(x,y)\hat{q}\right)_{y}\right)
dy\,dx\,dt = 0,
\end{align}
and if we can define an adjoint problem such that all but the first term 
in this equation vanishes then we are left with
\begin{align}\label{eq:reduced2d}
\int_a^b \int_{\alpha}^{\beta}  \hat{q}^T(x, y, t_f)q(x, y, t_f) dy\,dx
= \int_a^b \int_{\alpha}^{\beta} \hat{q}^T(x, y, t_0)q(x, y, t_0)dy\,dx,
\end{align}
which is the expression that allows us to use the inner product of the
adjoint and forward problems at each time step to determine what regions will
influence the point of interest at the final time. 

This requires that the adjoint equation have the form
\begin{equation}
 \hat{q}_{t} +
\left(A_1^T(x,y)\hat{q}\right)_{x} + \left(A_2^T(x,y)\hat{q}\right)_{y} = 0
\end{equation}
over the same domain as the forward problem. \revised{ and also
imposes coastline and boundary conditions similar to those that constrained 
the forward problem.}{
The correct boundary conditions to use are zero normal velocity at
all interfaces between any wet cell and dry cell so that the
boundary terms also drop out of expression \cref{eq:intbyparts2d} to obtain
\cref{eq:reduced2d}.}

\section{GeoClaw Tsunami Modeling Example}\label{sec:geoclawex}

\added{Finally, we present a tsunami propagation 
example utilizing the adjoint method to guide adaptive mesh refinement, as
implemented in the GeoClaw software package.
In principle the adjoint flagging methodology could be used in conjunction
with other tsunami models, although we know of no other open source software
that provides similar adaptive mesh refinement capabilities.}
Verification of the results
is performed by comparing the results from the adjoint method 
approach to the default approach already present in the GeoClaw software.
\added{The algorithms used in GeoClaw for tsunami modeling 
are described in detail in 
\cite{LeVequeGeorgeBerger2011}, and only a brief introduction will
be given here.
In general, GeoClaw solves the two-dimensional
nonlinear shallow water equations in the form of a nonlinear system of  
hyperbolic conservation laws for depth and momentum.
The numerical methods used are high-resolution Godunov-type finite
volume methods, in which the discrete solution is given by cell averages of
depth and momentum over the grid cells and the solution is updated in each
time step based on fluxes computed at the cell edges.  In Godunov-type
methods, these fluxes are determined by solving a ``Riemann problem'' at each cell
edge, which consists of the hyperbolic problem with piecewise constant initial
data given by the adjacent cell averages.  The general theory of these
methods is presented in \cite{LeVeque2002}, for example, and details of the
Riemann solver in GeoClaw can be found in \cite{dgeorge:jcp}.
Away from coastlines, this solver reduces to a Roe solver for
the shallow water equations plus bathymetry, which means that the
eigenstructure of a locally linearized Riemann problem is solved at each
cell interface, making it no more expensive in the deep ocean than simply
solving the linearized equations, but also capable of robustly 
handling nonlinearity near shore and inundation.}

As an example, we consider a tsunami generated 
by a hypothetical earthquake on the Alaska-Aleutian subduction zone, the
event denoted AASZe04 in the probabilistic hazard study of Crescent City, CA
performed by \cite{GonzalezLeVequeEtAl2014}. 
The major waves impinging on Crescent City from this 
hypothetical tsunami all occurred within 11 hours after the earthquake, 
so simulations will be run to this time.
To simulate the effects of this tsunami on Crescent City 
a coarse grid is used over the entire Pacific (1 degree resolution) where the ocean
is at rest. In addition to AMR being used to track propagating waves on finer grids, 
higher levels of refinement are allowed or enforced around Crescent City when the
tsunami arrives. A total of 4 levels of refinement are used, starting with 1-degree 
resolution on the coarsest level, and with refinement ratios of 5, 6, and 6 from one 
level to the next. Only 3 levels were allowed over most of the Pacific, and the remaining
level was used over the region around Crescent City. Level 4, with 20-second resolution,
is still too coarse to provide any real detail on the effect of the tsunami on the harbor. 
It does, however, allow for a comparison of flagging cells for refinement using the adjoint 
method and using the default method implemented in Geoclaw. 
\revised{}{In this simulation we used 1 arc-minute bathymetry from the
ETOPO1 Global Relief Model of \cite{etopo1}
for the entire simulation area, as well as 1 arc-second 
and 1/3 arc-second bathymetry over the region about Crescent City
from \cite{CCtopo}.
Internally, GeoClaw constructs a piecewise-bilinear function 
from the union of any provided topography files. This function 
is then integrated over computational grid cells to obtain a single
cell-averaged topography value in each grid cell in a manner that is
consistent between refinement levels.}

\begin{figure}[h!]
\begin{minipage}[c]{0.03\linewidth}
\hspace{0.1cm}
\end{minipage}
\begin{minipage}[c]{0.235\linewidth}
\centering
Grids
\end{minipage}
\begin{minipage}[c]{0.235\linewidth}
\centering
Solution
\end{minipage}
\begin{minipage}[c]{0.235\linewidth}
\centering
Grids
\end{minipage}
\begin{minipage}[c]{0.235\linewidth}
\centering
Solution
\end{minipage}\\
\begin{minipage}[c]{0.03\linewidth}
\begin{sideways}
t = 2 hours
\end{sideways}
\end{minipage}
\begin{minipage}[c]{0.235\linewidth}
\includegraphics[width=\textwidth]{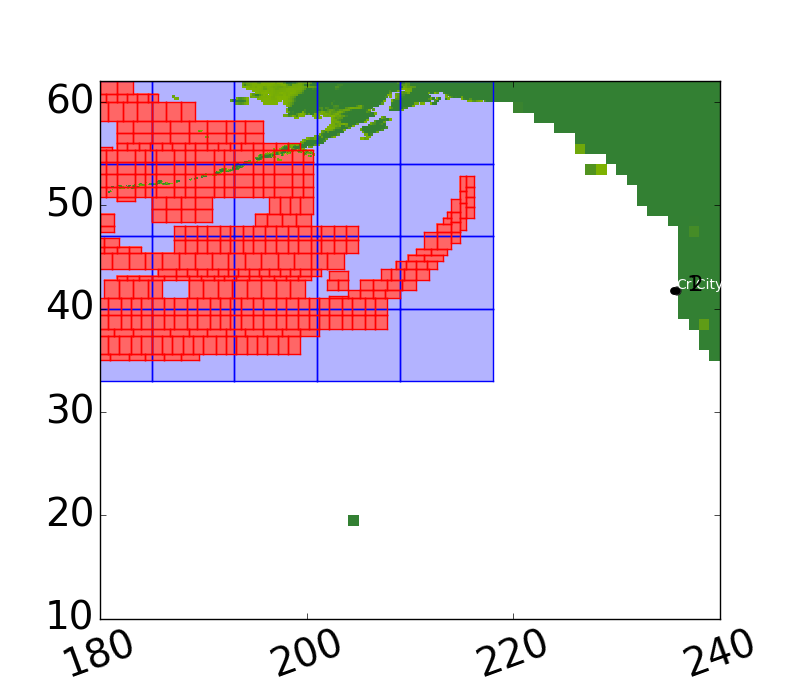}
\end{minipage}
\begin{minipage}[c]{0.235\linewidth}
\includegraphics[width=\textwidth]{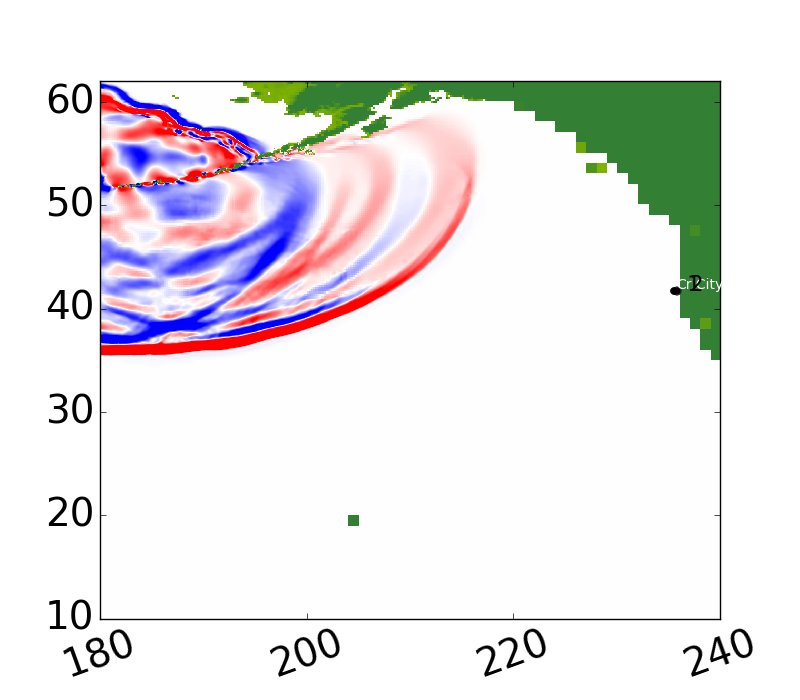}
\end{minipage}
\begin{minipage}[c]{0.235\linewidth}
\includegraphics[width=\textwidth]{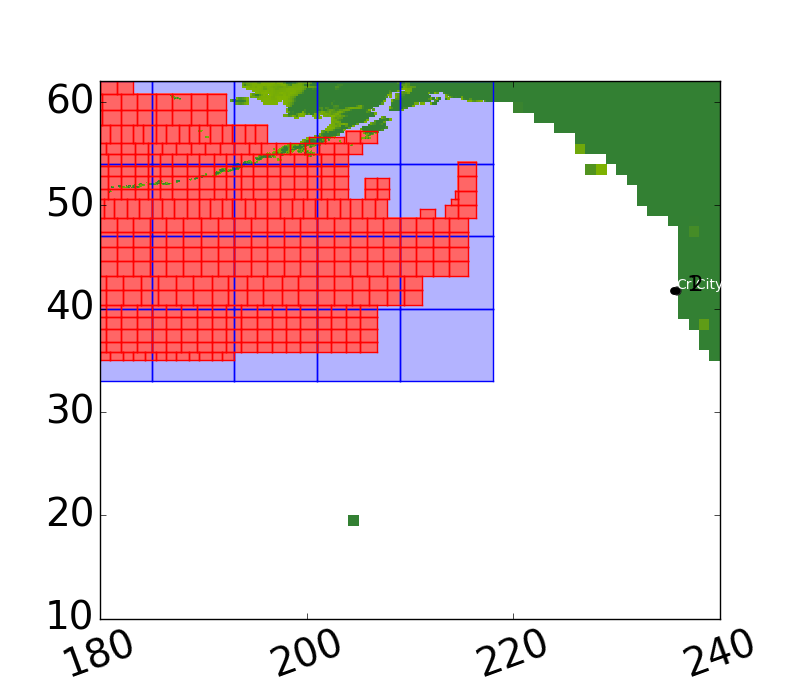}
\end{minipage}
\begin{minipage}[c]{0.235\linewidth}
\includegraphics[width=\textwidth]{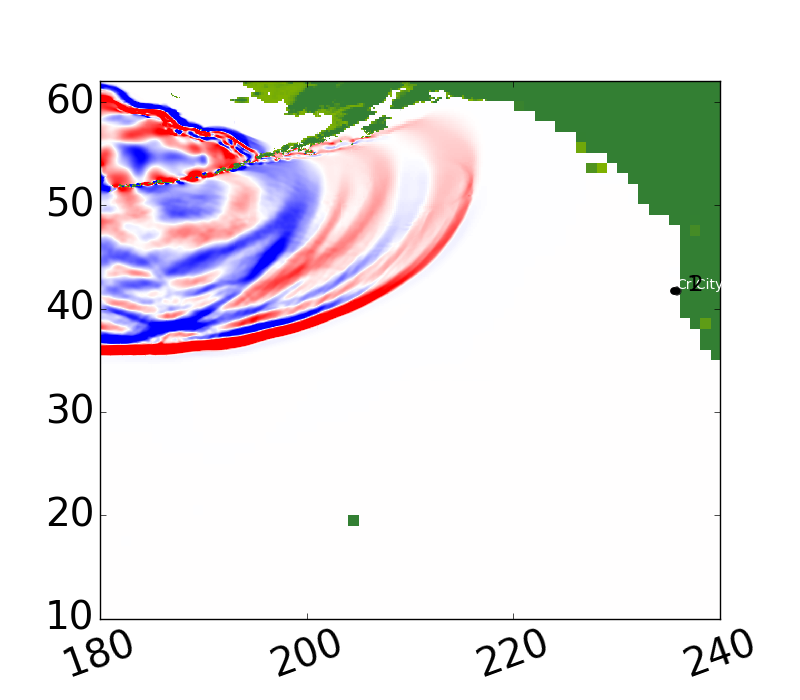}
\end{minipage}\\
\begin{minipage}[c]{0.03\linewidth}
\begin{sideways}
t = 4 hours
\end{sideways}
\end{minipage}
\begin{minipage}[c]{0.235\linewidth}
\includegraphics[width=\textwidth]{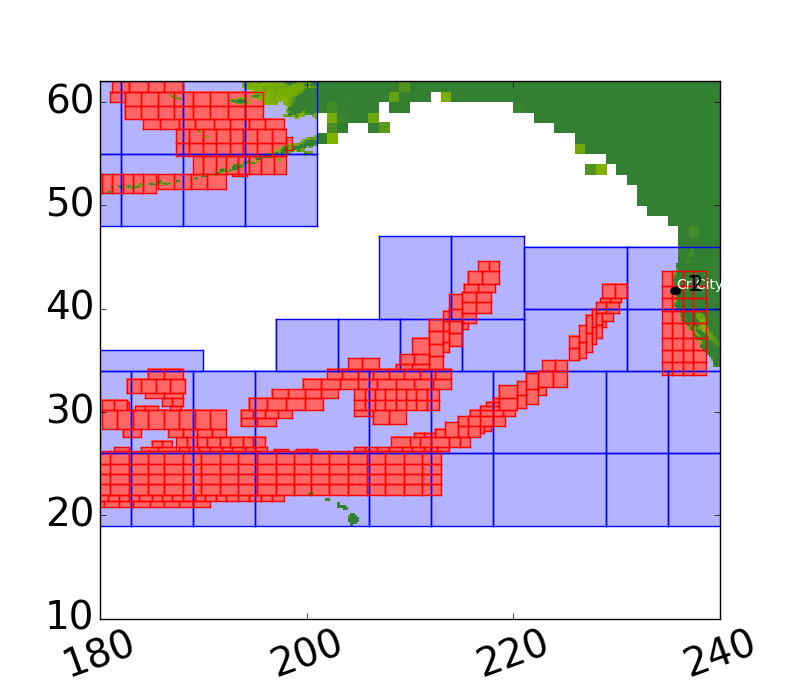}
\end{minipage}
\begin{minipage}[c]{0.235\linewidth}
\includegraphics[width=\textwidth]{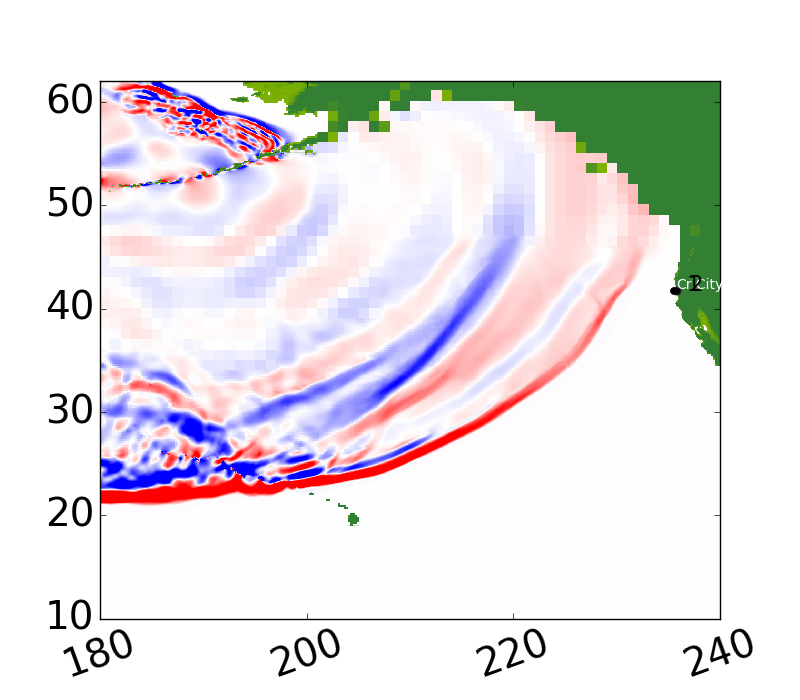}
\end{minipage}
\begin{minipage}[c]{0.235\linewidth}
\includegraphics[width=\textwidth]{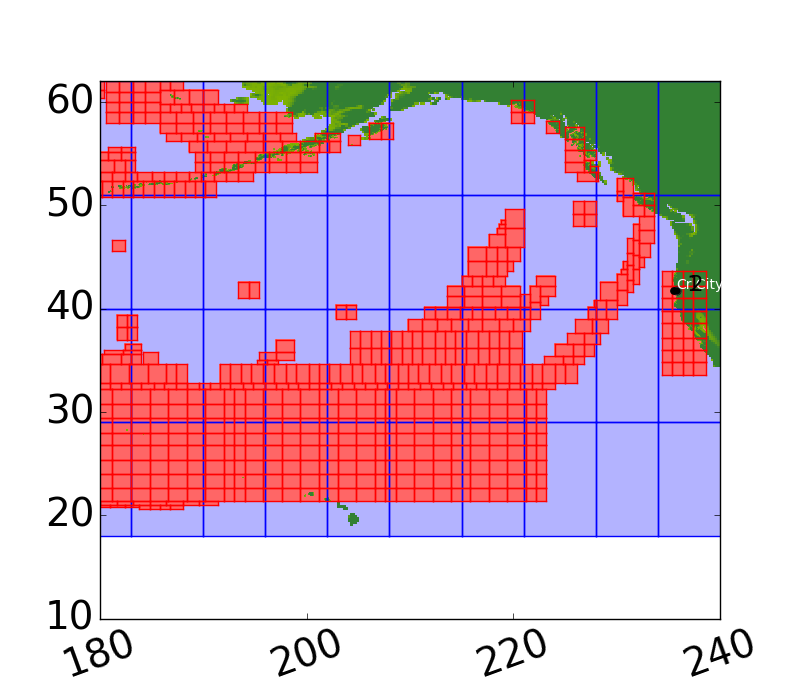}
\end{minipage}
\begin{minipage}[c]{0.235\linewidth}
\includegraphics[width=\textwidth]{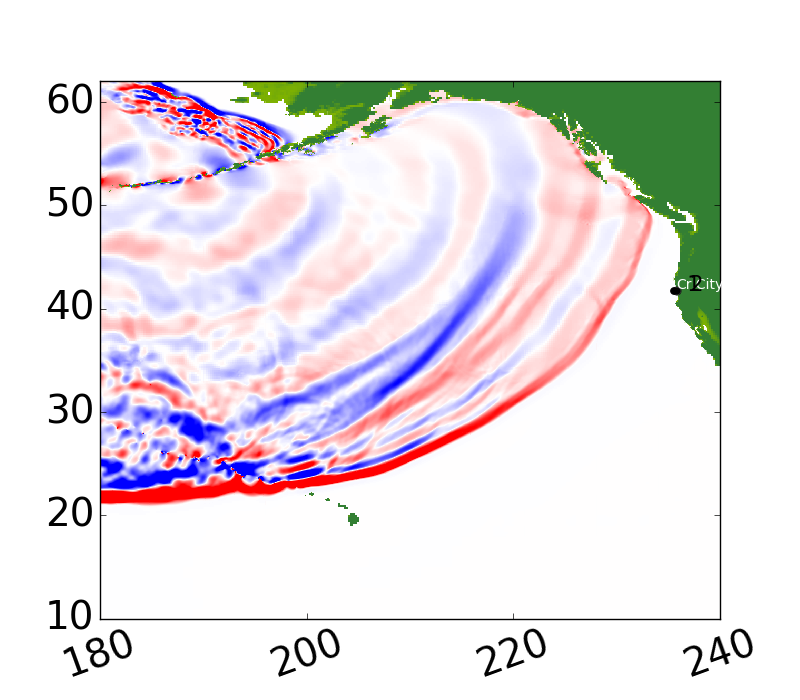}
\end{minipage}\\
\begin{minipage}[c]{0.03\linewidth}
\begin{sideways}
t = 6 hours
\end{sideways}
\end{minipage}
\begin{minipage}[c]{0.235\linewidth}
\includegraphics[width=\textwidth]{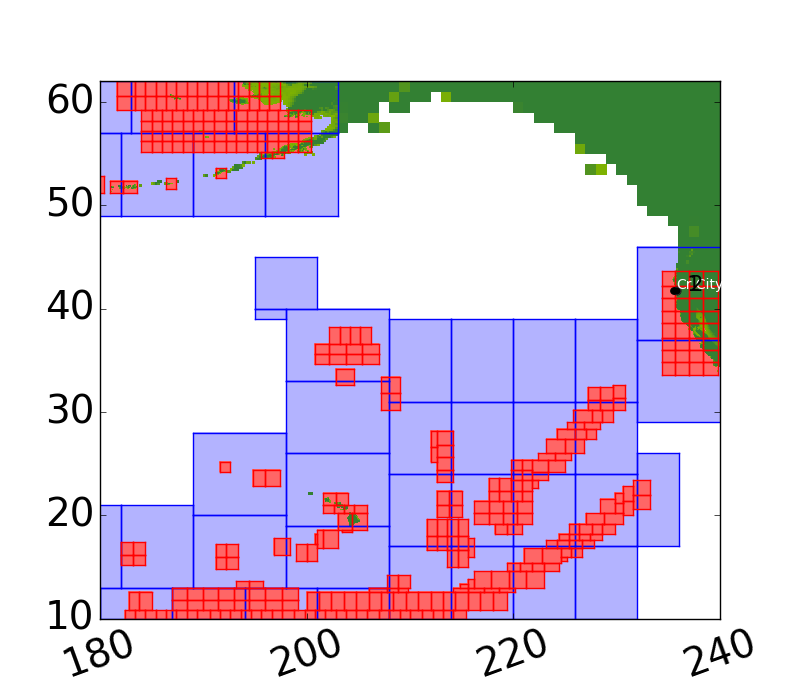}
\end{minipage}
\begin{minipage}[c]{0.235\linewidth}
\includegraphics[width=\textwidth]{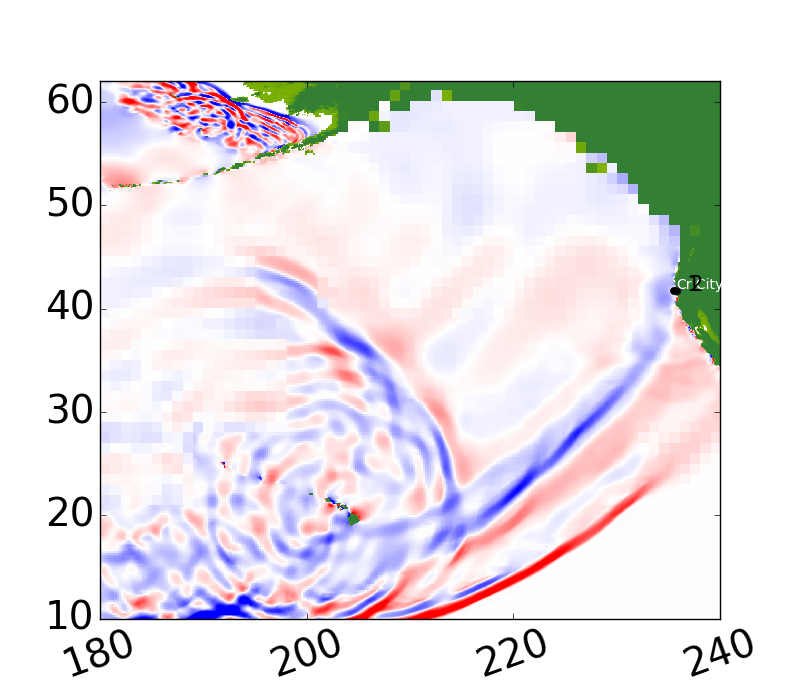}
\end{minipage}
\begin{minipage}[c]{0.235\linewidth}
\includegraphics[width=\textwidth]{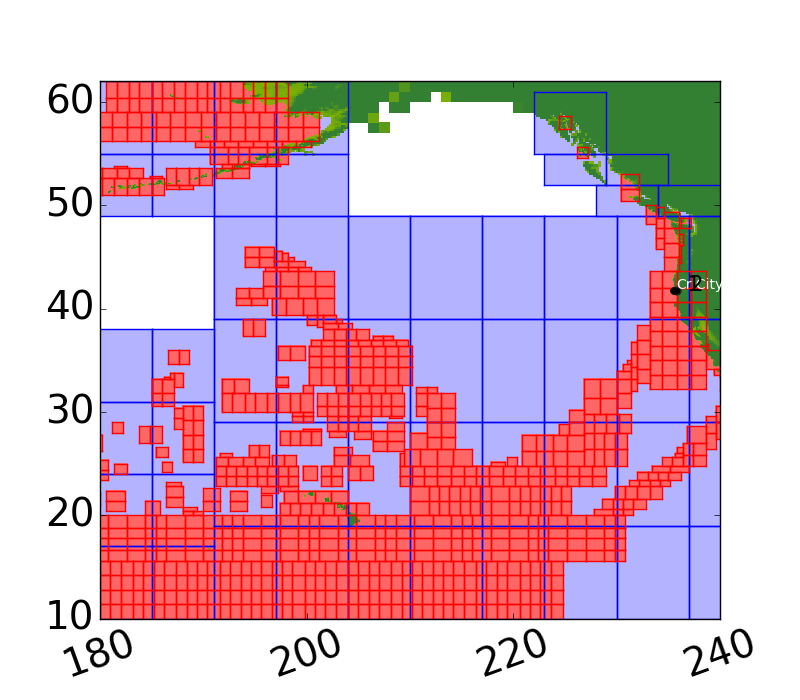}
\end{minipage}
\begin{minipage}[c]{0.235\linewidth}
\includegraphics[width=\textwidth]{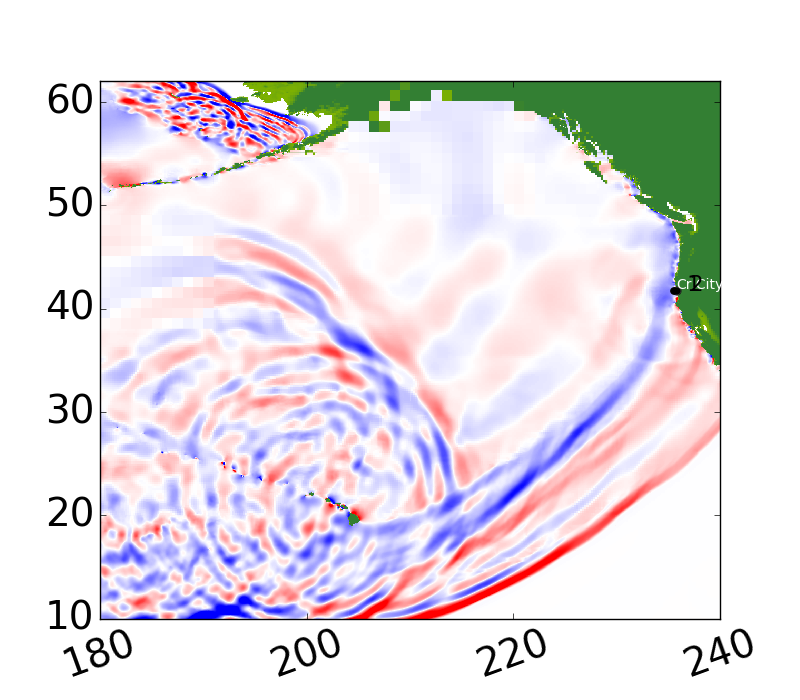}
\end{minipage}\\
\begin{minipage}[c]{0.03\linewidth}
\begin{sideways}
t = 8 hours
\end{sideways}
\end{minipage}
\begin{minipage}[c]{0.235\linewidth}
\includegraphics[width=\textwidth]{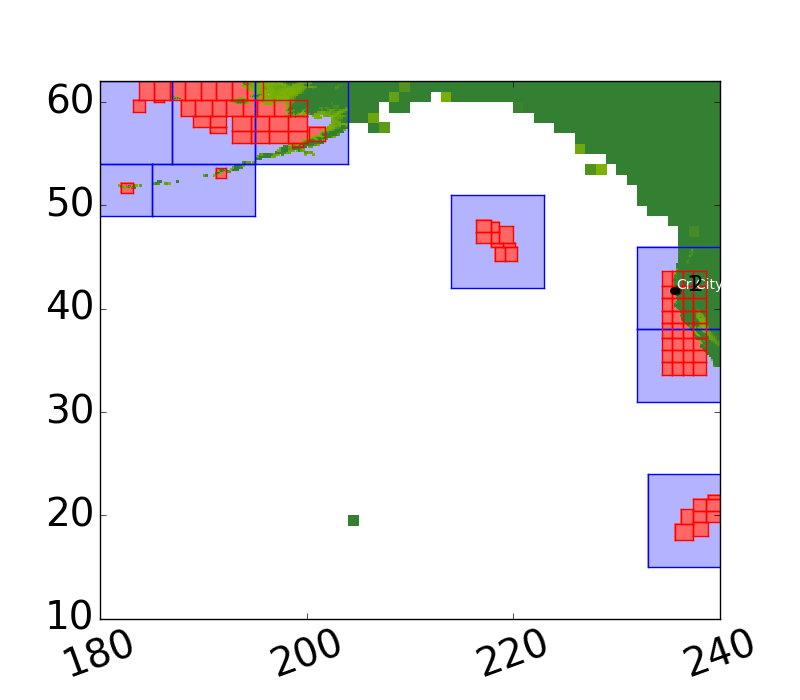}
\end{minipage}
\begin{minipage}[c]{0.235\linewidth}
\includegraphics[width=\textwidth]{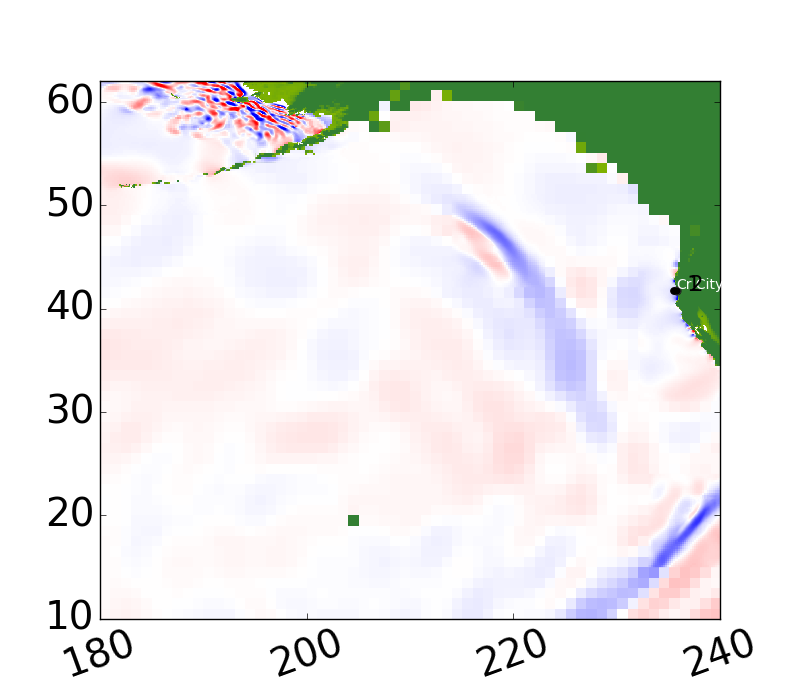}
\end{minipage}
\begin{minipage}[c]{0.235\linewidth}
\includegraphics[width=\textwidth]{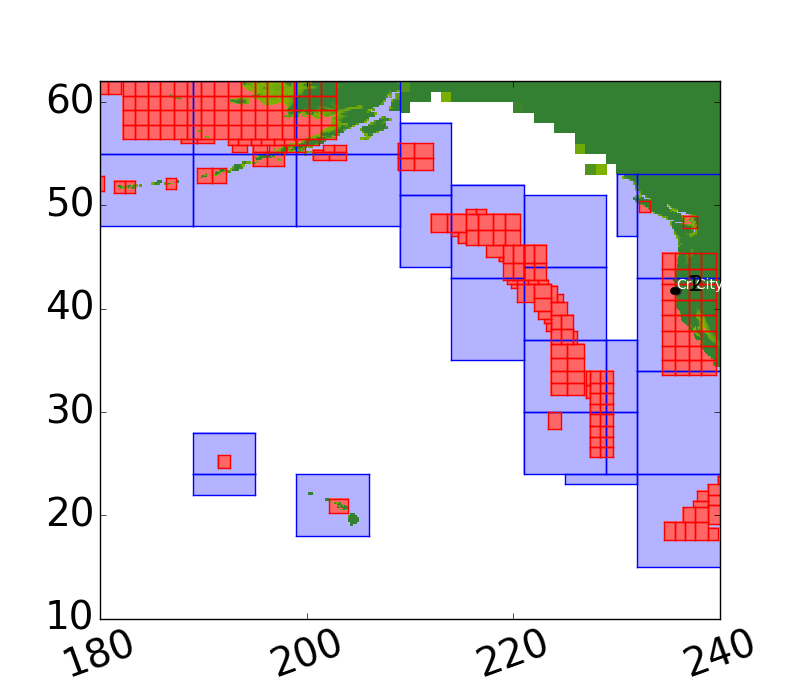}
\end{minipage}
\begin{minipage}[c]{0.235\linewidth}
\includegraphics[width=\textwidth]{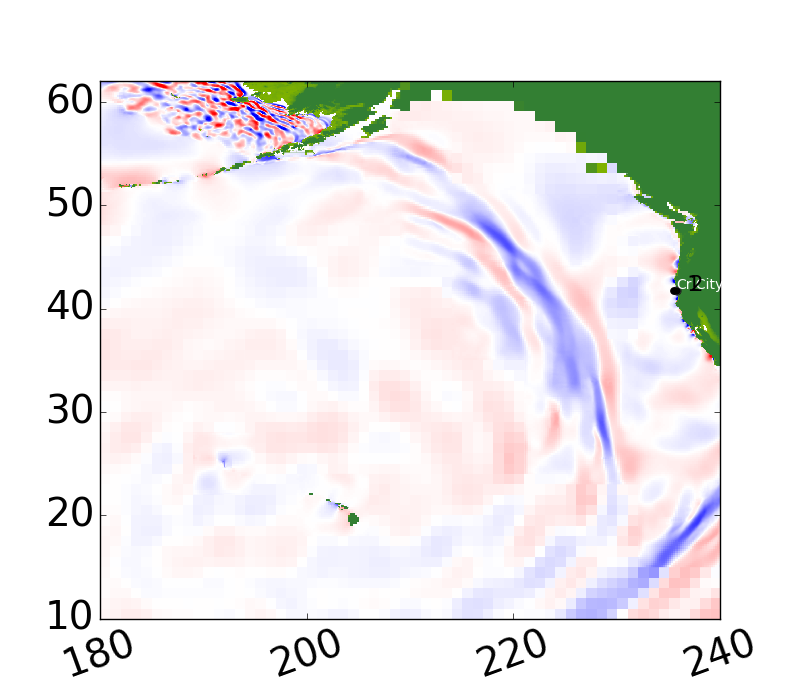}
\end{minipage}
 \caption{%
  Computed results for tsunami propagation problem on two different
  runs utilizing the surface-flagging technique. 
  The $x$-axis and $y$-axis are latitude and longitude, respectively. 
  The grids and solution along the left correspond to the simulation 
with a tolerance of $0.14$, and the \revised{girds}{grids} 
and solution along the right 
correspond to the simulation with a tolerance of $0.09$.
In the grid figures each color corresponds to a different level of refinement:
white for the coarsest level, blue for level two, and red for level three.
    The color scale for the solution figures goes from blue to red, and 
  ranges between $-0.3$ and 0.3 meters (surface elevation relative to sea
  level). Times are in hours after the earthquake.}
   \label{fig:Alaska_surfheight}
\end{figure}

The default Geoclaw refinement technique flags cells for refinement when the 
elevation of the sea surface relative to sea level is above some set tolerance, 
as described in \cite{LeVequeGeorgeBerger2011}, \added{where the adaptive
refinement and time stepping algorithms are described in more detail.} 
We will refer to this flagging method as 
surface-flagging. The value selected for this tolerance has a 
significant impact in the results calculated by the simulation, since a smaller tolerance
will result in more cells being flagged for refinement. Consequently, a smaller tolerance
both increases the simulation time required and theoretically increases the accuracy of
the results. 

Two Geoclaw simulations were \revised{preformed}{performed} using surface-flagging, 
one with a tolerance of $0.14$ and another with a tolerance of $0.09$. 
\Fig{fig:Alaska_surfheight} shows the results of these two simulations,
along with the \revised{girds}{grids} at different levels of refinement used, for the sake of comparison. 
Note that each grid outlined in the figure, colored based on its level of refinement, is 
a collection of cells at a particular refinement.
The grids and solution along the left of the figure correspond to the simulation 
with a tolerance of $0.14$, and the \revised{girds}{grids} and solution along the right of the figure 
correspond to the simulation with a tolerance of $0.09$.
 Note that 
the simulation with a surface-flagging tolerance of $0.14$ continues to refine the 
first wave until it arrives
at Crescent City about 5 hours after the earthquake, but after about 6 hours stops 
refining the main secondary wave which reflects off
the Northwestern Hawaiian (Leeward) Island chain before heading 
towards Crescent City. The second simulation, with a surface-flagging tolerance of 
$0.09$ continues to refine this secondary wave until it arrives at Crescent City. 

These two tolerances were selected because they are illustrative of two constraints 
that typically drive a Geoclaw simulation. The larger surface-flagging tolerance, 
of $0.14$ \revised{is}{was found to be} approximately the largest
tolerance that will refine the initial wave until it reaches Crescent City. Therefore, it 
essentially corresponds to a lower limit on the time required \added{by the
standard surface-flagging approach}: any Geoclaw simulation with a
larger surface-flagging tolerance would run more quickly but would fail to give accurate results for even
the first wave. Note that for this particular example, even when the simulation will only
give accurate results for the first wave to reach Crescent City, a large area of the wave 
front that is not headed directly towards Crescent City is being refined with the AMR. 
The smaller surface-flagging tolerance, of $0.09$, refines all of the waves of interest that
impinge on Crescent City thereby giving more accurate results at the expense
of longer computational time, 
\added{particularly since it also refines waves that will never reach Crescent City.}

Now we consider the adjoint approach, which will allow us to refine only those
sections of the wave that will affect our target region.
For this example, we are interested in the accurate calculation of the 
water surface height in the area about Crescent City, California. To focus on this
area, we define a circle of radius $1^\circ$
centered about $(x_c,y_c) = (235.80917,41.74111)$ where $x$ and $y$ are being
measured in degrees. Setting
\begin{align*}
J = \int_{x_{min}}^{x_{max}}\int_{y_{min}(x)}^{y_{max}(x)}\eta(x,y,t_f)dy\,dx,
\end{align*}
where the limits of integration define the appropriate circle,
the problem then requires that
\begin{align}\label{eqn:phi_Alaska}
\varphi (x,y) = \left[ \begin{matrix}
I (x,y) \\ 0 \\ 0
\end{matrix}
\right], 
\end{align}
where
\begin{align}\label{eqn:delta_Alaska}
I (x,y) = \left\{
     \begin{array}{ll}
       1 & \hspace{0.3in}\textnormal{if~} \sqrt{\left( x - x_c\right)^2
        + \left( y - y_c\right)^2} \leq 1,\\
       0 & \hspace{0.3in}\textnormal{otherwise.}
     \end{array}
   \right. 
\end{align}

\begin{figure}[!htbp]
 \centering
  \subfigure[$t_f -$ 1 hour]{
  \includegraphics[width=0.225\textwidth]{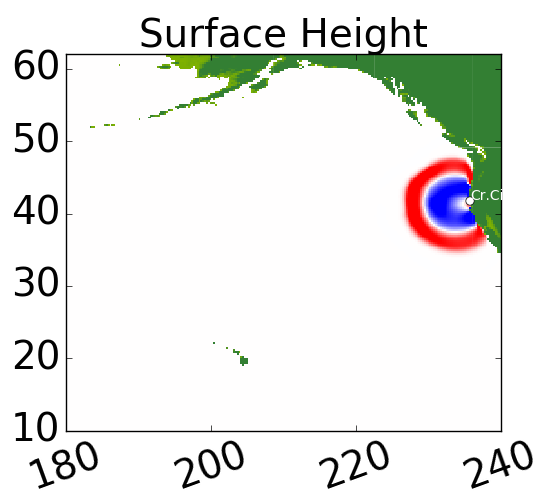}
   \label{fig:Alaska_a1}
   }
   \subfigure[$t_f -$ 3 hours]{
  \includegraphics[width=0.225\textwidth]{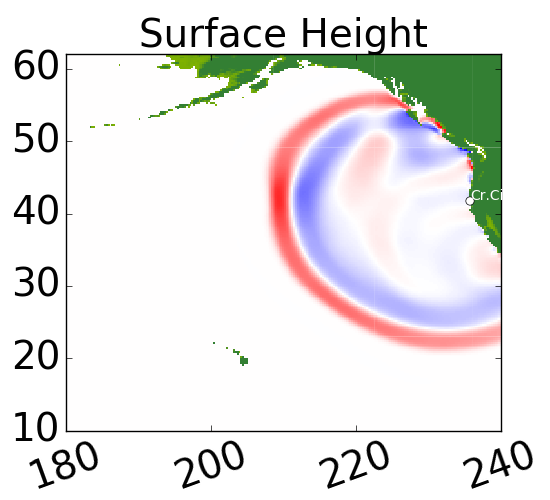}
   \label{fig:Alaska_a2}
   }
  \subfigure[$t_f -$ 5 hours]{
  \includegraphics[width=0.225\textwidth]{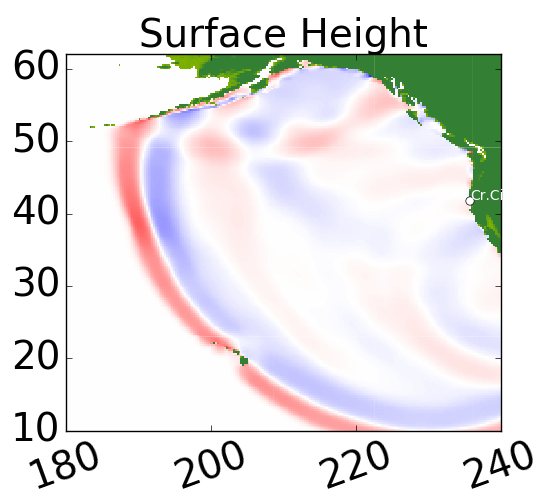}
   \label{fig:Alaska_a3}
   }
  \subfigure[$t_f -$ 7 hours]{
  \includegraphics[width=0.225\textwidth]{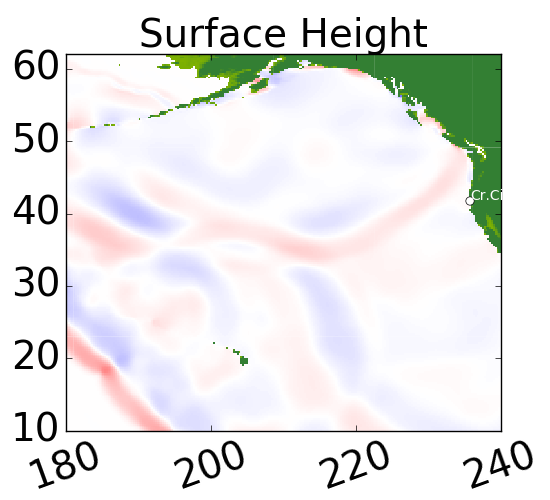}
   \label{fig:Alaska_a4}
   }
 \caption{%
   Computed results for tsunami propagation adjoint problem. Times shown are 
   the number of hours before the final time, since the ``initial'' conditions are 
   given at the final time. 
  The color scale goes from blue to red, and ranges between $-0.05$ and $0.05$.}
   \label{fig:Alaska_adjoint}
\end{figure}

This function $\varphi(x,y)$ defines the ``initial data'' $\hat{q}(x,y,t_f)$
for the adjoint problem.
\Fig{fig:Alaska_adjoint} shows the results for the 
simulation of this adjoint problem. For this simulation a grid with 15
arcminute $= 0.25^\circ$
resolution was used over the entire Pacific and no grid refinement was allowed. 
The simulation was run out to 11 hours.

The topography files used for the adjoint problem are the same as those used
for the forward problem.
However, given that the adjoint problem is being solved on a coarser grid than 
the forward problem, the coastline between the two simulations varies.
Since the coastline varies between the 
two simulations, when computing the inner product it is possible to find grid cells
that are wet in the forward solution and dry in the adjoint solution. In this case, 
the inner product in those grid cells is set to zero.
 
The simulation of this tsunami using adjoint-flagging for the AMR
 was run using the same initial grid over the Pacific, the same refinement 
ratios, and the same initial water displacement as our previous 
surface-flagging simulations. The only difference between this simulation
and the previous one is the flagging technique utilized. 
The first waves arrive at Crescent
City around 4 hours after the earthquake, so we set $t_s = 3.5$ hours and $t_f = 11$ 
hours. Recall that the areas where the maximum magnitude of the 
inner product over the appropriate time range,
\begin{align*}
\max\limits_{T \leq \tau \leq t}
\left|\hat{q}^T(x,y,\tau ) q(x,y,t)\right|
\end{align*}
with $T = \min (t + t_f - t_s, 0)$, is large are the 
areas where adaptive mesh refinement should take place.

\begin{figure}[!htbp]
\begin{minipage}[c]{0.02\linewidth}
\hspace{0.1cm}
\end{minipage}
\begin{minipage}[c]{0.32\linewidth}
\centering
Grids
\end{minipage}
\begin{minipage}[c]{0.32\linewidth}
\centering
Solution
\end{minipage}
\begin{minipage}[c]{0.32\linewidth}
\centering
Inner Product
\end{minipage}\\
\begin{minipage}[c]{0.02\linewidth}
\begin{sideways}
t = 2 hours
\end{sideways}
\end{minipage}
\begin{minipage}[c]{0.32\linewidth}
\includegraphics[width=\textwidth]{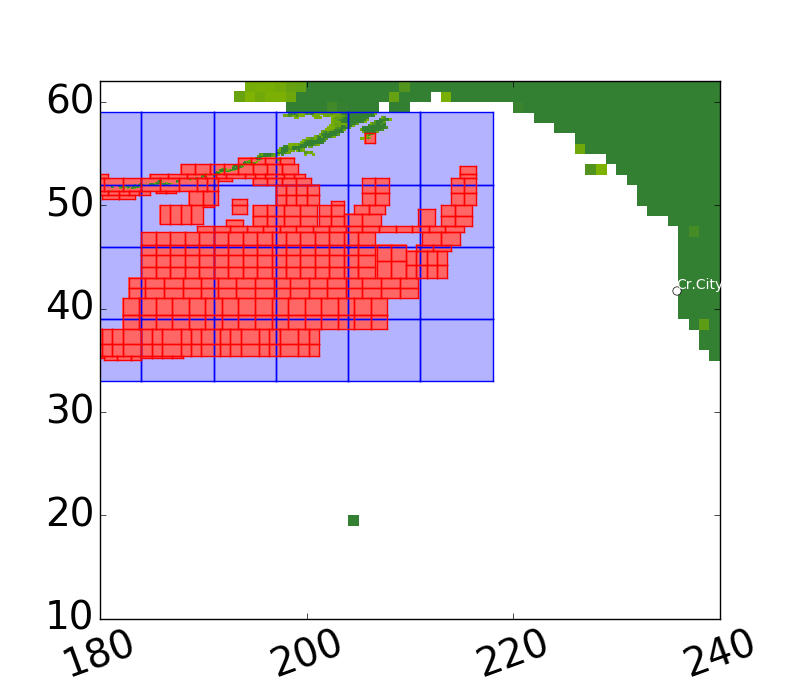}
\end{minipage}
\begin{minipage}[c]{0.32\linewidth}
\includegraphics[width=\textwidth]{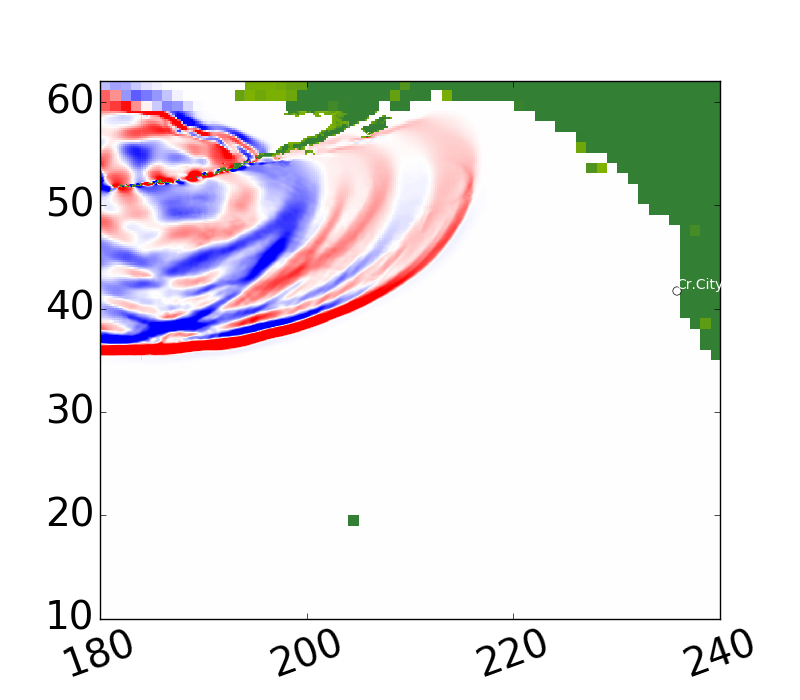}
\end{minipage}
\begin{minipage}[c]{0.32\linewidth}
\includegraphics[width=\textwidth]{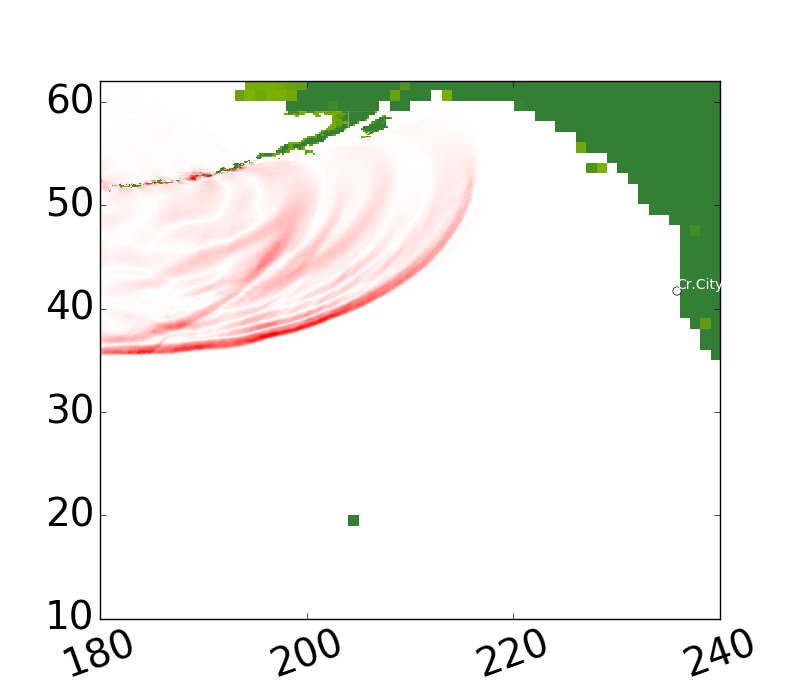}
\end{minipage}\\
\begin{minipage}[c]{0.02\linewidth}
\begin{sideways}
t = 4 hours
\end{sideways}
\end{minipage}
\begin{minipage}[c]{0.32\linewidth}
\includegraphics[width=\textwidth]{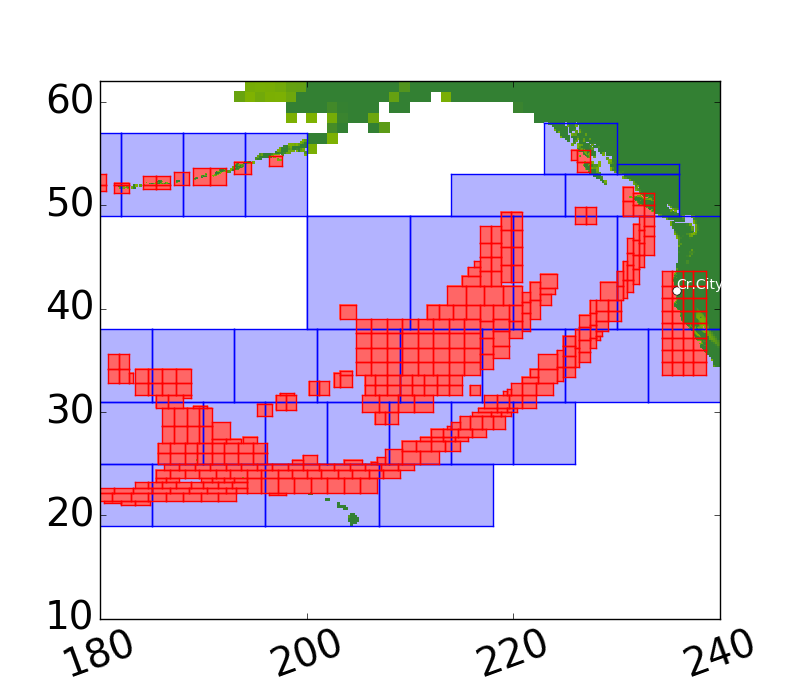}
\end{minipage}
\begin{minipage}[c]{0.32\linewidth}
\includegraphics[width=\textwidth]{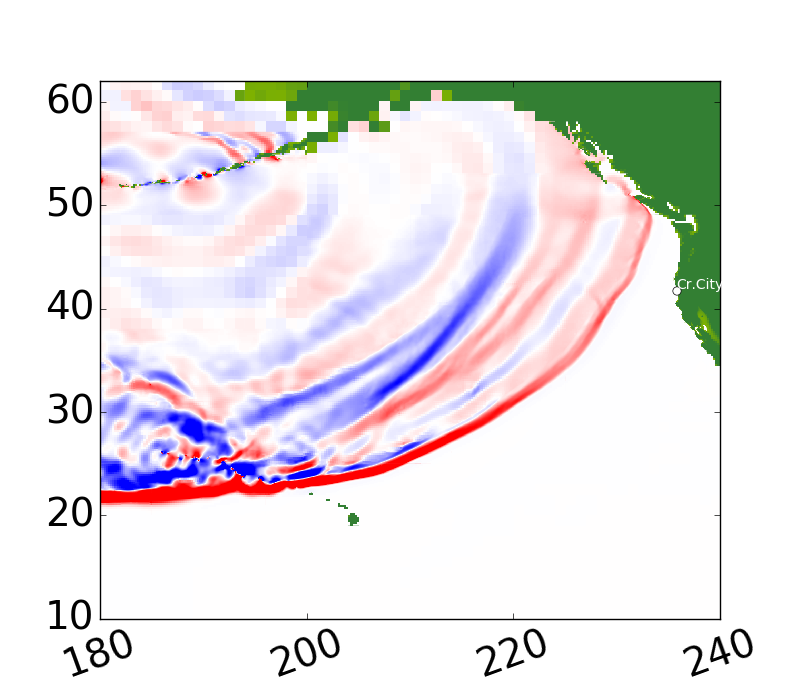}
\end{minipage}
\begin{minipage}[c]{0.32\linewidth}
\includegraphics[width=\textwidth]{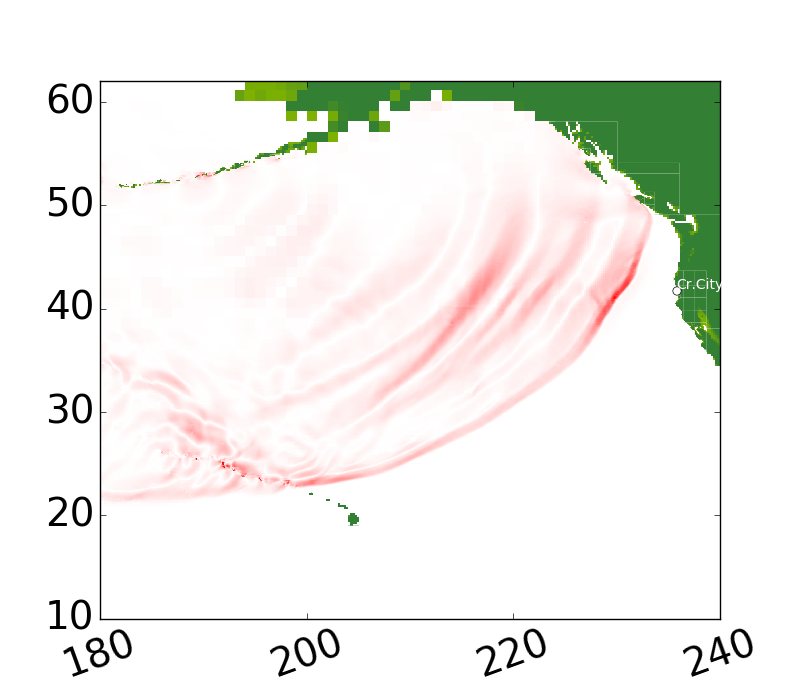}
\end{minipage}\\
\begin{minipage}[c]{0.02\linewidth}
\begin{sideways}
t = 6 hours
\end{sideways}
\end{minipage}
\begin{minipage}[c]{0.32\linewidth}
\includegraphics[width=\textwidth]{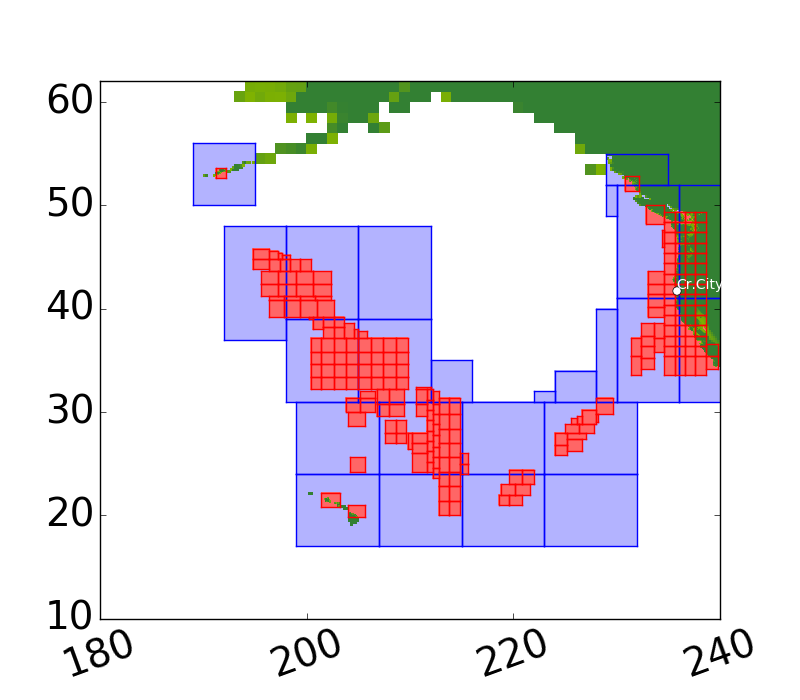}
\end{minipage}
\begin{minipage}[c]{0.32\linewidth}
\includegraphics[width=\textwidth]{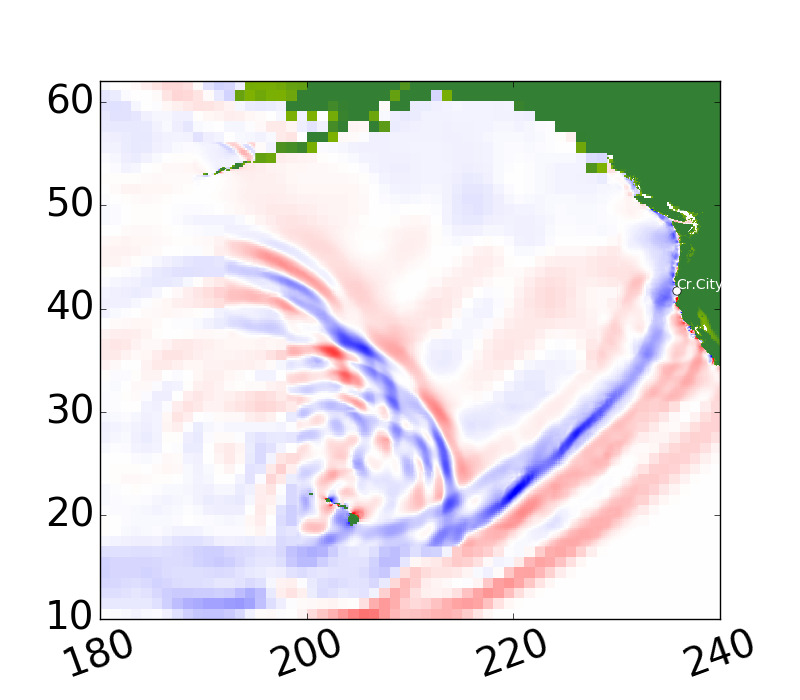}
\end{minipage}
\begin{minipage}[c]{0.32\linewidth}
\includegraphics[width=\textwidth]{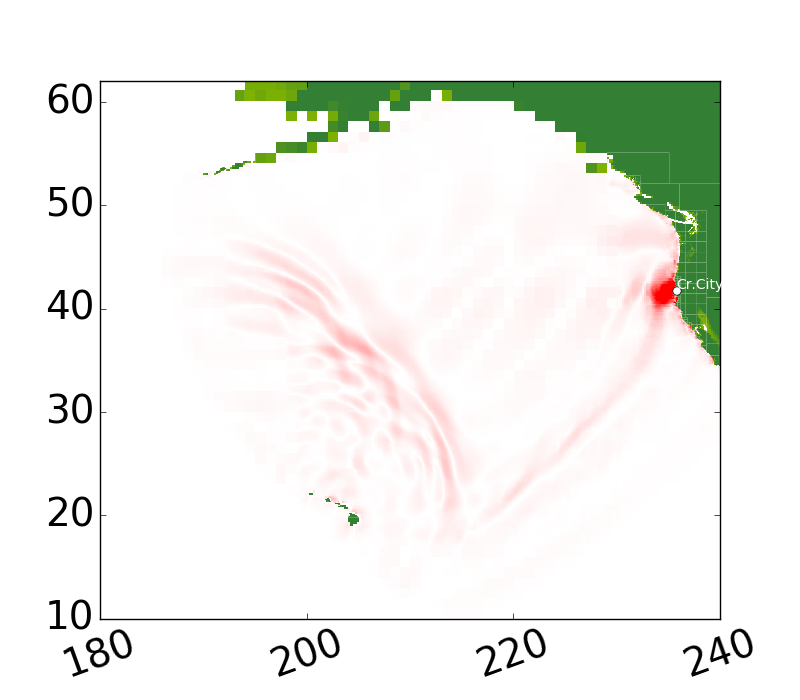}
\end{minipage}\\
\begin{minipage}[c]{0.02\linewidth}
\begin{sideways}
t = 8 hours
\end{sideways}
\end{minipage}
\begin{minipage}[c]{0.32\linewidth}
\includegraphics[width=\textwidth]{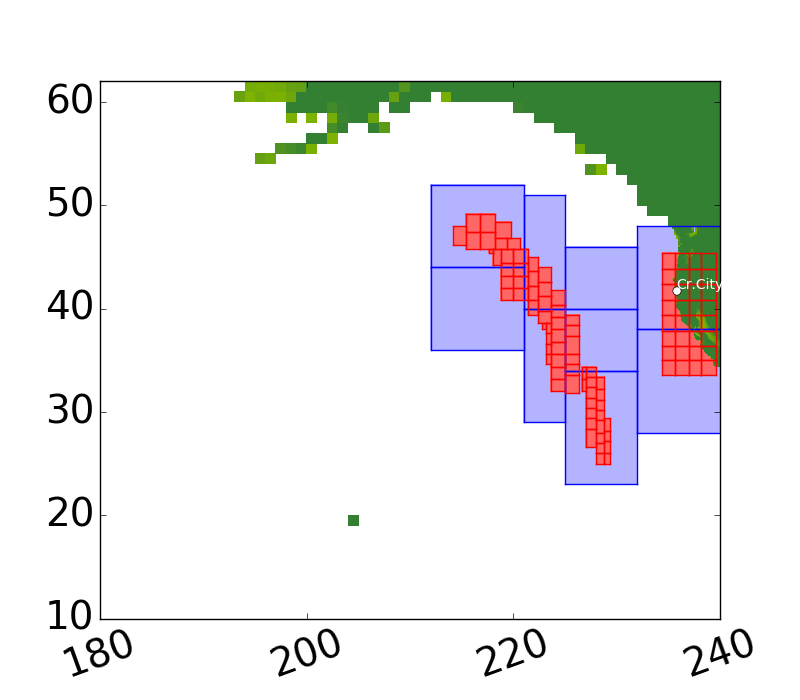}
\end{minipage}
\begin{minipage}[c]{0.32\linewidth}
\includegraphics[width=\textwidth]{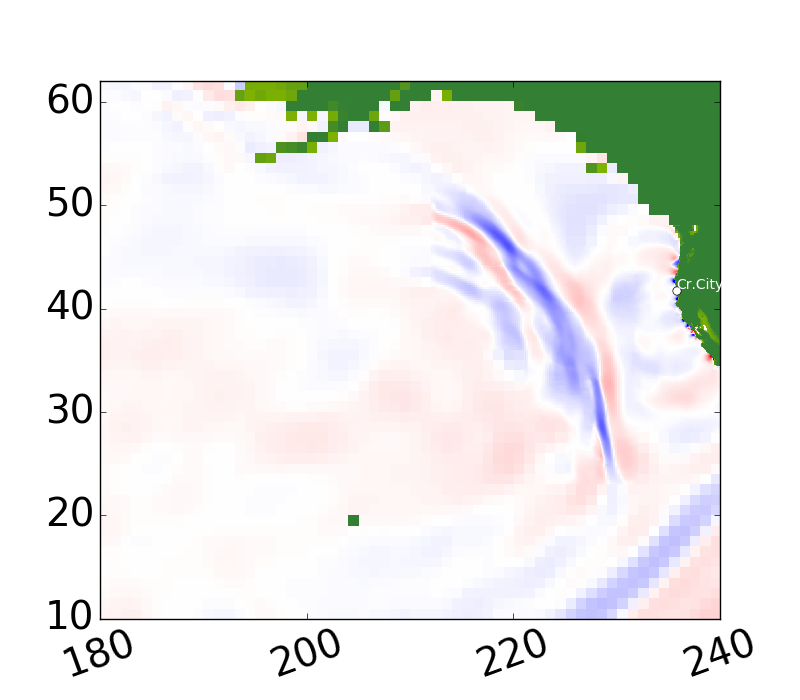}
\end{minipage}
\begin{minipage}[c]{0.32\linewidth}
\includegraphics[width=\textwidth]{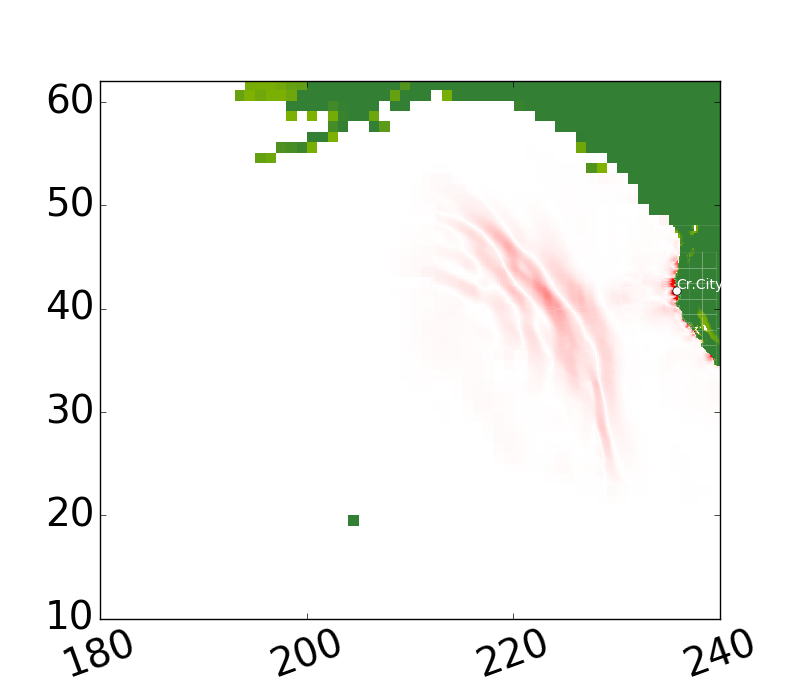}
\end{minipage}
 \caption{%
  Computed results for tsunami propagation problem when adjoint-flagging is
  used. 
  The $x$-axis and $y$-axis are latitude and longitude, respectively. 
  In the grid figures each color corresponds to a different level of refinement:
white for the coarsest level, blue for level two, and red for level three.
    The color scale for the surface height figures goes from blue to red, and 
  ranges between $-0.3$ and 0.3. The color scale for the inner 
  product figures goes from white to red, and ranges between 
  0 and 0.04.
  Times are in hours since the earthquake.}
   \label{fig:Alaska}
\end{figure}

\Fig{fig:Alaska} shows the Geoclaw results for the surface height at various different times, 
along with the grids at different levels of refinement that were used and
 the maximum inner product in the appropriate time range. 
 Compare this figure to \Fig{fig:Alaska_surfheight}, noting the extent of the grids 
 at each refinement level for each of the three simulations.

\subsection{Computational Performance}\label{sec:tsunamiperf}
The above example was run on a quad-core laptop, for both the surface-flagging 
and adjoint-flagging methods, and the OpenMP option of GeoClaw was enabled 
which allowed all four cores to be utilized.
The timing results for the tsunami simulations are shown in Table \ref{tb:timing_tsunami}.
 Recall that two simulations were run using 
 surface-flagging, one with a tolerance of $0.14$ (``Large Tolerance'' 
 in the table) and another with a tolerance of $0.09$ (``Small Tolerance'' in the table).
  Finally, a GeoClaw example using adjoint-flagging was run with a 
  tolerance of $0.004$ (``Forward'' in the table), which of course required 
  a simulation of the adjoint problem the timing for which is also shown in the table.

  As expected, between the two GeoClaw simulations which utilized  
surface-flagging the one with the larger tolerance took significantly less time. 
Note that although 
  solving the problem using adjoint-flagging did require two different simulations, 
  the adjoint problem and the forward problem, the computational time required is
only slightly more than the timing required for the large tolerance
surface-flagging case.
 \begin{table}[htbp]
\caption{Timing comparison for the example in Section \ref{sec:tsunamiperf} given in seconds. }
\begin{center}\footnotesize
\renewcommand{\arraystretch}{1.3}
\begin{tabular}{| c  | c | c | c|}\hline
\multicolumn{2}{ |c| }{Surface-Flagging} &\multicolumn{2}{ |c| }{Adjoint-Flagging}\\\hline
 \bf Small Tolerance & \bf Large Tolerance &\bf Forward & \bf Adjoint\\ \hline
 8310 &  5724 & 5984 & 27 \\\hline
\end{tabular}
\end{center}
\label{tb:timing_tsunami}
\end{table}
Another consideration when comparing the adjoint-flagging method with the
surface-flagging method already in place in GeoClaw is the accuracy of the
results. To test this, gauges were placed in the example and the output
at the gauges compared across the two different methods. 

\begin{figure}[h!]
 \centering
  \subfigure{
  \includegraphics[width=0.9\textwidth]{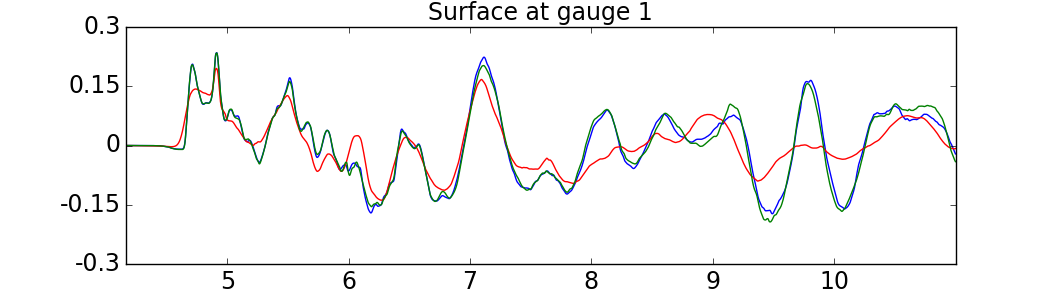}
   \label{fig:Alaska_gauge1}
   }
   \subfigure{
  \includegraphics[width=0.9\textwidth]{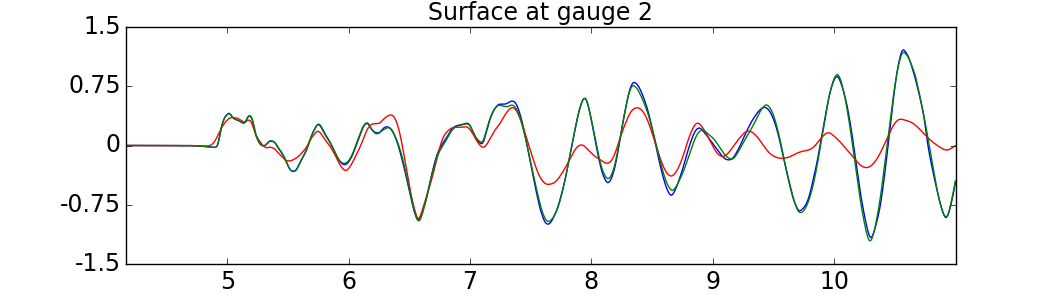}
   \label{fig:Alaska_gauge2}
   }
 \caption{%
  Computed results at gauges for tsunami propagation problem.
  The results from the simulation using the adjoint method are shown in 
  blue, the results from the simulation using the surface-flagging method with 
  a tolerance of $0.14$ are shown in red, and the results from the simulation 
  using the surface-flagging method with a tolerance of $0.09$ are shown in green.
  Along the $x$-axis, the time since the occurrence of the earthquake is shown in hours. 
}
   \label{fig:Alaska_gauges}
\end{figure}

For the tsunami example two gauges are used: gauge 1 is placed at 
$(x,y) = (235.536, 41.67)$ which is on the continental shelf to the west of Crescent 
City, and gauge 2 is placed at $(x,y) = (235.80917,41.74111)$ which is in the 
harbor of Crescent City. In \Fig{fig:Alaska_gauges} the gauge results from the adjoint
 method are shown in blue, the results from the surface-flagging technique with a 
tolerance of $0.14$ are shown in red, and the results from the surface-flagging 
technique with a tolerance of $0.09$ are shown in green. Note that the blue and green
lines are in fairly good agreement, indicating that the 
use of the adjoint method achieved a comparable accuracy 
with the smaller tolerance run
using the surface-flagging method, although the time required was 
significantly less.
While the larger tolerance run using the surface-flagging method had a 
similar time requirement to the adjoint method simulation, it agrees 
fairly well only for the first wave but then rapidly loses accuracy. 

\section{Conclusions}\label{sec:conclusions}

\added{In this paper, we first presented the adjoint methodology in some detail in hopes
that it will also be useful in other tsunami modeling software in the
future, and perhaps in other contexts for exploring sensitivities.}

Integrating the adjoint method approach to cell flagging into the already existing AMR 
algorithm in GeoClaw results in significant time savings for the tsunami simulation shown 
in this work. 
For the tsunami example we examined here, the surface-flagging simulation with a tolerance of
$0.14$ has the advantage of a low computational time but only provides accurate results 
for the first wave to reach Crescent City. The surface-flagging simulation
with a tolerance of $0.09$ refines more waves and therefore provides accurate results for 
a longer period of time. However, it has the disadvantage of having a long computational 
time requirement \added{since it refines in many regions where the waves are
not important to the modeling of Crescent City}. 
The use of the adjoint method allows us to retain accurate results 
while also reducing the computational time required. \added{Equally, or
perhaps even more importantly, the adjoint approach gives 
more confidence that the appropriate waves have been
refined to capture the tsunami impact at the target location than {\em ad
hoc} attempts to guide the refinement regions manually.}
 
\revised{
 Therefore, using the adjoint method to guide adaptive mesh refinement can reduce the 
 computational expense of solving a system of equations while retaining the accuracy
 of the results by enabling targeted 
 refinement of the regions of the domain that will influence a specific area
of interest.}{}

The code for all the examples presented in this work is available online 
at \cite{adjointCode}, and includes
the code for solving the adjoint Riemann problems. 
This code can be easily modified to solve
\revised{different}{other tsunami modeling} problems. 
\revised{Note that it requires the addition of a new Riemann solver
for the adjoint problem.}{This repository also contains other examples
illustrating how adjoint flagging can be used with AMRClaw, the more general
adaptive refinement code in Clawpack for general hyperbolic systems.
Another example of the adjoint method being used in GeoClaw for tsunami 
modeling can be found in \cite{Borrero2015}.}
 
\section{Additional Comments and Future Work}\label{sec:futurework}
The method described in this paper flags cells for refinement wherever
the magnitude of the relevant inner product between the
forward and adjoint solutions is above some tolerance.  Choosing a
sufficiently small tolerance will trigger refinement of
all regions where the forward solution might need to be refined, and in our
current approach these will be refined to the finest level specified in the
computation \added{for refinement in the ocean (finer levels may be imposed
near the target location).}  
\added{We believe this approach is already
a significant advance over the method currently used in GeoClaw for many
trans-oceanic modeling problems.  Although a tolerance must be chosen in
order to define the cutoff for flagging cells based on the inner product,
this is similar to the current need for setting a tolerance on the surface
elevation for flagging cells, and has the great advantage that it identifies
the waves that will reach the target location rather than potentially
refining everywhere there are waves.}

\revised{To}{We will also continue to explore ways to better}
optimize efficiency --- it would be desirable to have error bounds based on the
adjoint solution that could be used to refine the grid more selectively to
achieve some target error tolerance for the final quantity of interest.
\revised{We}{For some problems we} believe this can be accomplished using  
the Richardson extrapolation error estimator that is built into 
AMRClaw to estimate the point-wise error in the forward solution
and then using the adjoint solution to estimate its effect on the final 
quantity of interest. 
This is currently under investigation and we hope to develop 
a robust strategy that can be applied to a wide variety of problems
for general inclusion into Clawpack.
\added{However, there are potential difficulties in deriving more precise error
estimation via the adjoint method for use in GeoClaw.  The functional $J(q)$
that we have defined, e.g. by integrating the solution against
the piecewise constant function over a one degree
square around Crescent City defined by \cref{eqn:delta_Alaska}, 
is not exactly the quantity we are trying to compute in the end.  Rather we
wish to compute the time history at one or more particular gauge locations
(as shown in \Fig{fig:Alaska_gauges}) or the detailed inundation in the
community.  The functional $J(q)$ is simply designed to radiate waves from
the Crescent City vicinity in order to determine what waves in the forward
problem are important to track.  The error in some true quantity of interest
such as the maximum flooding depth may vary from point to point in the 
community.  For this reason more research is required to investigate the
extent to which error estimates on $J(q)$ can be employed in practice.
Another potential difficulty is that the bathymetry $B(x,y)$ is not at all
smooth at the grid resolution typically needed to model trans-ocean wave
propagation. Because the wave length of tsunamis is so long, good accuracy
is often observed in spite of this (as found in many validation studies of
GeoClaw and other tsunami software based on the shallow water equations).
However, this may limit the applicability of Richardson extrapolation error
estimation.}

\added{Using adaptive refinement in solving the adjoint equation may also be
desirable.}
In this paper the adjoint solution was computed on a fixed 
 grid. Allowing AMR to take place when solving the adjoint equation
\revised{should}{could} increase the accuracy 
 of the results, since it would enable a more accurate \added{(while still
efficient)} evaluation of the inner product between the 
 forward and adjoint solutions. 
 In an effort to guide the AMR of the adjoint problem in a similar manner
 to the method used for the forward problem, the two problems would need to be solved somewhat 
 in conjunction and the inner product between the two considered for both the flagging of the 
 cells in the adjoint problem as well as the flagging of cells in the forward problem.
 One approach that is used to tackle this issue is checkpointing, where the forward
 problem is solved and the solution at 
 a small number of time steps is stored for use when solving the 
 adjoint problem, as seen for example in \cite{WangMoinIaccarino2009}.
 The automation of this process is another area of future work, 
 and involves developing an evaluation 
 technique for determining the number of checkpoints to save or 
 when to shift from refining the adjoint solution to refining the forward 
 solution and vice versa. 
 
For the tsunami example presented here
we linearized the shallow water equations about the ocean at rest, 
and the adjoint equations then have essentially the same form.
This is sufficient for many important applications, in
particular for tsunami applications where the goal is to track waves in the
ocean where the \revised{linearized equtions are very
accurate}{linearization is essentially independent of the forward solution.}
If an adjoint
method is desired in the inundation zone, or for other nonlinear hyperbolic
equations,  then the 
adjoint equation is derived by linearizing about
a particular forward solution.  This would 
again require the development of an automated process 
to shift between solving the forward problem, linearizing about that forward problem, and solving 
the corresponding adjoint problem.
Finally, in this work we assumed wall boundary conditions when a wave interacted
with the coastline in the adjoint problem. This assumption, along with the use
 of the linearized shallow water equations, becomes significant when a wave 
 approaches a shore line. Allowing for more accurate \revised{iterations}{interactions}
  between waves 
 and the coastline in the solution of the adjoint problem is another area for 
 future work.


\end{document}